\pdfoutput=1
\documentclass[11pt,a4paper,reqno]{amsart}

\usepackage[british]{babel}

\usepackage[DIV=9,oneside,BCOR=0mm]{typearea}
\usepackage{microtype}
\usepackage[T1]{fontenc}
\usepackage{setspace}
\usepackage{amssymb}
\usepackage{amsmath}
\usepackage{amsthm}
\usepackage{enumitem}
\usepackage{mathrsfs}
\usepackage[colorlinks,citecolor=blue,urlcolor=black,linkcolor=blue,linktocpage]{hyperref}
\usepackage{mathtools}
\usepackage{xcolor}
\usepackage{xfrac}

\usepackage{accents}

\newcommand{\dtilde}[1]{\tilde{\raisebox{0pt}[0.85\height]{$\tilde{#1}$}}}

\newtheorem{innercustomgeneric}{\customgenericname}
\providecommand{\customgenericname}{}
\newcommand{\newcustomtheorem}[2]{%
	\newenvironment{#1}[1]
	{%
		\renewcommand\customgenericname{#2}%
		\renewcommand\theinnercustomgeneric{##1}%
		\innercustomgeneric
	}
	{\endinnercustomgeneric}
}

\newcustomtheorem{customthm}{Theorem}
\newcustomtheorem{customcor}{Corollary}

\newtheorem{thm}{Theorem}[section]

\newtheorem{cor}[thm]{Corollary}
\newtheorem{lem}[thm]{Lemma}
\newtheorem{prop}[thm]{Proposition}

\theoremstyle{definition}

\newtheorem{definition}[thm]{Definition}

\newtheorem{remark}[thm]{Remark}

\renewcommand{\epsilon}{\varepsilon}
\renewcommand{\phi}{\varphi}
\newcommand{\defeq}{\mathrel{\mathop:}=}

\DeclareMathOperator{\Hom}{Hom}
\DeclareMathOperator{\Iso}{Iso}
\DeclareMathOperator{\width}{width}
\DeclareMathOperator{\Nest}{\mathcal{N}}
\DeclareMathOperator{\Nestmax}{\mathcal{N}_{max}}
\DeclareMathOperator{\alg}{V\!}
\DeclareMathOperator{\rk}{rk}
\DeclareMathOperator{\rank}{rank}
\DeclareMathOperator{\ind}{ind}
\DeclareMathOperator{\M}{M}
\DeclareMathOperator{\I}{I}
\DeclareMathOperator{\lat}{L}
\DeclareMathOperator{\E}{E}
\DeclareMathOperator{\GL}{GL}
\DeclareMathOperator{\PGL}{PGL}
\DeclareMathOperator{\SL}{SL}
\DeclareMathOperator{\cent}{Z}
\DeclareMathOperator{\Cl}{Cl}
\DeclareMathOperator{\id}{id}
\DeclareMathOperator{\im}{im}

\DeclareMathOperator{\A}{A}

\DeclareMathOperator{\R}{\mathbb{R}}

\DeclareMathOperator{\N}{\mathbb{N}}

\DeclareMathOperator{\ball}{B}
\DeclareMathOperator{\trap}{T}
\DeclareMathOperator{\diam}{diam}
\DeclareMathOperator{\Neigh}{\mathcal{U}}
\DeclareMathOperator{\comp}{\mathbf{c}}

\DeclareMathOperator{\Pfin}{\mathcal{P}_{fin}}

\begin{document}

\setlist{noitemsep}

\author{Friedrich Martin Schneider}
\address{F.M.~Schneider, Institute of Discrete Mathematics and Algebra, TU Bergakademie Freiberg, 09596 Freiberg, Germany}
\email{martin.schneider@math.tu-freiberg.de}
\thanks{This research is funded by the Deutsche Forschungsgemeinschaft (DFG, German Research Foundation) -- Projektnummer 561178190}

\title[Geometric properties of unit groups of continuous rings]{Geometric properties of unit groups of \\von Neumann's continuous rings}
\date{\today}

\begin{abstract} 
We prove that, if $R$ is a non-discrete irreducible, continuous ring, then its unit group $\GL(R)$, equipped with the topology generated by the rank metric, is topologically simple modulo its center, path-connected, locally path-connected, bounded in the sense of Bourbaki, and not admitting any non-zero escape function. All these topological insights are consequences of more refined geometric results concerning the rank metric, in particular with regard to the set of algebraic elements. Thanks to the phenomenon of automatic continuity, our results also have non-trivial ramifications for the underlying abstract groups.
\end{abstract}



\maketitle

\allowdisplaybreaks


\tableofcontents

\newpage

\section{Introduction}

The study of continuous rings, initiated by John von Neumann in his discovery of a continuous-dimensional analogue of finite-dimensional projective geometry~\cite{VonNeumannBook}, has gained new momentum in recent years motivated by a variety of questions pertaining to group rings~\cite{ElekSzabo,linnell,elek2013}, operator algebras~\cite{LinnellSchick,elek}, and topological dynamics~\cite{SchneiderGAFA,SchneiderIMRN,BernardSchneider}. The present manuscript aims to advance this development by establishing three results of geometric flavor about the unit groups of such rings.

A \emph{continuous ring} is a regular ring $R$ such that the corresponding lattice of its principal right ideals constitutes a continuous geometry (see Section~\ref{section:continuous.rings} for details). Among a number of other profound results, von Neumann established in~\cite{VonNeumannBook} that every irreducible, continuous ring $R$ admits a unique \emph{rank function} $\rk_{R} \colon R \to [0,1]$. The associated \emph{rank metric} $d_{R} \colon R \times R \to [0,1], \, (a,b) \mapsto \rk_{R}(a-b)$ furnishes $R$ with a compatible topology---the \emph{rank topology}. Moreover, as the rank metric restricts to a bi-invariant metric on the unit group $\GL(R)$, i.e., the group of multiplicatively invertible elements of $R$, the pair $(\GL(R),{d_{R}})$ is an instance of a \emph{metric group} (see Remark~\ref{remark:properties.pseudo.rank.function}\ref{remark:unit.group}), around whose geometric properties this paper revolves. A common theme in our three main results lies in the application of the following theorem by von Neumann~\cite{VonNeumann37,Halperin62}: if $R$ is non-discrete with respect to the rank topology, then the set $\A(R)$ of elements of $R$ \emph{algebraic} over its center $\cent(R)$ is dense in $R$. This marvelous insight turns out to be particularly powerful when being used in combination with the \emph{continuous triangularization theorem} of~\cite[Theorem~9.11]{SchneiderGAFA}.

Our first result concerns the problem of densely covering a topological group $G$ by products of conjugacy classes, i.e., of subsets of the form $\Cl_{G}(g) \defeq \left. \! \left\{ hgh^{-1} \, \right\vert h \in G \right\}$ where $g \in G$. In particular, given a single (non-central) element $g \in G$, one may ask whether $G = \overline{\Cl_{G}(g)^{n}}$ for some natural number $n$. For unit groups of non-discrete irreducible, continuous rings, we answer this question by establishing a concrete bound depending on an element's rank distance from the center.

\begin{customthm}{A}[Theorem~\ref{theorem:tbng}]\label{theorem:a} Let $R$ be a non-discrete irreducible, continuous ring. If $m \in \N$ and $g_{1},\ldots,g_{m} \in \GL(R)$ are such that $\sum\nolimits_{i=1}^{m} d_{R}(g_{i},\cent (R)) > 12$, then \begin{displaymath}
	\GL(R) \, = \, \overline{\Cl_{\GL(R)}(g_{1})\cdots \Cl_{\GL(R)}(g_{m})} .
\end{displaymath} In particular, if $g \in \GL(R)\setminus \cent (R)$ and $m \in \N$ are such that $m > \tfrac{12}{d_{R}(g,\cent (R))}$, then \begin{displaymath}
	\GL(R) \, = \, \overline{\Cl_{\GL(R)}(g)^{m}} .
\end{displaymath} \end{customthm}

For the \emph{projective unit group} of a non-discrete irreducible, continuous ring $R$, i.e., the quotient topological group \begin{displaymath}
	\PGL(R) \, \defeq \, \GL(R)/\cent(\GL(R)) \, \stackrel{\ref{proposition:center}}{=} \, \GL(R)/(\cent(R)\setminus \{ 0 \}) ,
\end{displaymath} Theorem~\ref{theorem:a} entails the following.

\begin{customcor}{A}[Corollary~\ref{corollary:pgl}]\label{corollary:a} Let $R$ be a non-discrete irreducible, continuous ring. Then $\PGL(R)$ has topological bounded normal generation. In particular, $\PGL(R)$ is topologically simple. \end{customcor}

For the particular example of $R$ being the completion of the inductive limit of the diagonally embedded matrix rings $\M_{2^{n}}(F)$ $(n \in \N)$ with respect to the joint extension of the normalized rank metrics over a finite field $F$ (see Section~\ref{section:escape.dynamics}), this corollary had been established by Carderi and Thom in~\cite{CarderiThom} using work of Liebeck and Shalev~\cite{LiebeckShalev}. For metric ultraproducts of irreducible, continuous rings, Theorem~\ref{theorem:a} admits a proper strengthening (Corollary~\ref{corollary:ultraproduct}), which generalizes a result by Gismatullin, Majcher and Ziegler~\cite[Theorem~9.4]{GMZ25} concerning metric ultraproducts of special linear groups over the field of complex numbers. Via different methods, the author's more recent work~\cite[Corollary~1.5]{SchneiderSimple} establishes bounded normal generation (in particular, simplicity) of the group $\PGL(R)$ for every non-discrete irreducible, continuous ring, thus strengthening Corollary~\ref{corollary:a} and answering two questions by Carderi and Thom~\cite[p.~260]{CarderiThom}.

The second major topic of this work is a construction of geodesics resulting from continuous triangularization of algebraic elements. For a metric group $(G,d)$, let us define $\Delta(G,d)$ to be the set of elements of $G$ connected to the identity by a geodesic.

\begin{customthm}{B}[Theorem~\ref{theorem:geodesic}]\label{theorem:b} Let $R$ be a non-discrete irreducible, continuous ring. Then \begin{displaymath}
	\GL(R) \cap \A(R) \, \subseteq \, \Delta(\GL(R),{d_{R}}) .
\end{displaymath} In particular, \begin{displaymath}
	\GL(R) \, = \, \overline{\Delta(\GL(R),{d_{R}})} .
\end{displaymath} \end{customthm}

This insight comes about as a consequence of a more comprehensive study of geodesics in algebraic sets (Lemma~\ref{lemma:geodesic}). Our techniques have further applications: for instance, it turns out that, for any algebraic subset $X \subseteq R$ intersecting the center of~$R$, the metric space $(X,d_{R})$ is star-shaped (Proposition~\ref{proposition:star}). Among the direct ramifications of Theorem~\ref{theorem:b}, we deduce the following.

\begin{customcor}{B}[Corollary~\ref{corollary:length.space}]\label{corollary:b} Let $R$ be a non-discrete irreducible, continuous ring. Then the metric space $(\GL(R),d_{R})$ is a length space. In particular, $\GL(R)$ is both path-connected and locally path-connected with respect to the rank topology. \end{customcor}

Moreover, our work leads to new \emph{width bounds} for the considered unit groups. The \emph{width}~\cite{Liebeck} of a group $G$ with respect to a subset $S \subseteq G$ is defined as \begin{displaymath}
	\width(G,S) \, \defeq \, \inf \{ n \in \N \mid G = S^{n} \} \,  \in \, {\N} \cup {\{ \infty \}} .
\end{displaymath} Width questions have been attracting considerable attention for a long time, for example concerning the sets of involutions, commutators, and general verbal expressions (see, e.g.,~\cite{SegalBook} for more details). We define the \emph{geodesic width} of a metric group $(G,d)$ as the width of $G$ with respect to $\Delta(G,d)$. Combining a local version of Theorem~\ref{theorem:b}, namely Lemma~\ref{lemma:locally.algebraic}, with a recent decomposition theorem~\cite[Theorem~1.1]{BernardSchneider}, we establish a uniform upper bound on the geodesic width of the unit group of an arbitrary non-discrete irreducible, continuous ring (Theorem~\ref{theorem:geodesic.width}).

Our third main objective is a geometric observation concerning the escape dynamics, i.e., the behavior of \emph{escape functions} (Definition~\ref{definition:escape.function}), on unit groups of non-discrete irreducible, continuous rings. The study of such functions on general topological groups was initiated in~\cite{SchneiderSolecki24}, deeply inspired by Tao's account~\cite{TaoBook} on the solution to Hilbert's fifth problem. For a subset $U$ of a group $G$, we let $\trap(U)$ denote the union of all subgroups of $G$ contained in $U$. Moreover, the closed ball of radius $r \in \R_{>0}$ centered at the neutral element of a metric group $(G,d)$ will be denoted by $\ball_{r}(G,d)$.

\begin{customthm}{C}[Theorem~\ref{theorem:bounded}]\label{theorem:c} Let $R$ be a non-discrete irreducible, continuous ring, and let $n \in \N_{>0}$. Then \begin{displaymath}
	\GL(R) \cap \A(R) \, \subseteq \, \trap\!\left(\ball_{1/n}(\GL(R),d_{R})\right)^{n} .
\end{displaymath} In particular, \begin{displaymath}
	\GL(R) \, = \, \overline{\trap\!\left(\ball_{1/n}(\GL(R),d_{R})\right)^{n}} \, = \, \trap\!\left(\ball_{1/n}(\GL(R),d_{R})\right)^{n+1} .
\end{displaymath} \end{customthm}
	
Our Theorem~\ref{theorem:c}	implies that the topological unit group of a non-discrete irreducible, continuous ring does not possess any non-zero escape function (Corollary~\ref{corollary:no.escape}) and thus does not admit a non-trivial continuous homomorphism into any Hausdorff topological group with the escape property (Corollary~\ref{corollary:homomorphism.rigidity}). The latter excludes a broad range of target groups: as established in~\cite[Proposition~3.6]{SchneiderSolecki24}, the class of topological groups with the escape property includes all Banach--Lie groups, all locally compact groups, all non-archimedean topological groups, and all isometry groups of locally compact separable metric spaces. Combining this with the automatic continuity result of~\cite[Theorem~1.3]{BernardSchneider}, we moreover deduce the following consequences for the underlying abstract groups (not carrying any topology).

\begin{customcor}{C}[Corollary~\ref{corollary:rigidity}]\label{corollary:c} Let $R$ be a non-discrete irreducible, continuous ring. \begin{enumerate}
	\item Every action of $\GL(R)$ by isometries on a separable locally compact metric space is trivial. In particular, any action of $\GL(R)$ on a countable set is~trivial.
	\item Every homomorphism from $\GL(R)$ to any $\sigma$-compact locally compact topological group is trivial. In particular, $\GL(R)$ is minimally almost periodic.
\end{enumerate} \end{customcor}

Furthermore, the conclusion of Theorem~\ref{theorem:c} entails boundedness in the sense of Bourbaki (Corollary~\ref{corollary:bounded}), which allows us to improve upon a result by Carderi and Thom~\cite{CarderiThom} regarding triviality of representations on Hilbert spaces for unit groups of certain natural examples of continuous rings (Corollary~\ref{corollary:strongly.exotic}).

This manuscript is organized as follows. We start off with two preliminary sections, providing some background from metric geometry in Section~\ref{section:metric.geometry} and collecting the necessary prerequisites on continuous rings in Section~\ref{section:continuous.rings}. In Section~\ref{section:topological.simplicity}, we establish Theorem~\ref{theorem:a} along with its reinforcement for metric ultraproducts. The subsequent Section~\ref{section:geodesic.width} is dedicated to a study of rank-metric geodesics in continuous rings and their relation with the underlying algebraic geometry, which results in the proof of Theorem~\ref{theorem:b} and a discussion of various topological and algebraic ramifications. In the final Section~\ref{section:escape.dynamics} we turn to escape dynamics, proving Theorem~\ref{theorem:c} along with the above-mentioned corollaries.

\section{Some metric geometry}\label{section:metric.geometry}

The purpose of this section is to provide some basic background material from metric geometry, as it may be found, for instance, in~\cite{BridsonHaefliger,BuragoBuragoIvanov}.

If $(X,d)$ is a pseudo-metric space, then we let \begin{displaymath}
	d(x,Y) \, \defeq \, \inf \{ d(x,y) \mid y \in Y \}
\end{displaymath} for any $x \in X$ and $Y \subseteq X$. Now, let $(X,d)$ be a metric space. We say that $(X,d)$ has \begin{itemize}
	\item[---\,] \emph{midpoints} if, for all $x,y \in X$, there is $z \in X$ with $d(x,z) = \tfrac{1}{2}d(x,y) = d(z,y)$,
	\item[---\,] \emph{approximate midpoints} if, for all $x,y \in X$ and $\epsilon \in \R_{>0}$, there exists $z \in X$ such that $\left\lvert d(x,z) - \tfrac{1}{2}d(x,y) \right\rvert \leq \epsilon$ and $\left\lvert d(z,y) - \tfrac{1}{2}d(x,y) \right\rvert \leq \epsilon$.
\end{itemize} The \emph{length} of a path $\gamma$ in $(X,d)$, i.e., a continuous map $\gamma \colon [0,t] \to X$ with $t \in \R_{\geq 0}$, is defined as \begin{displaymath}
	\left. \ell_{d}(\gamma) \, \defeq \, \sup \! \left\{ \sum\nolimits_{i=0}^{n-1} d(\gamma(t_{i}),\gamma(t_{i+1})) \, \right\vert n \in \N_{>0}, \, 0 = t_{0} \leq \ldots \leq t_{n} = t \right\} .
\end{displaymath} A \emph{geodesic} in $(X,d)$ from one point $x \in X$ to another point $y \in Y$ is a $d$-isometric map $\phi \colon [0,d(x,y)] \to X$ with $\phi(0) = x$ and $\phi(d(x,y)) = y$. The metric space $(X,d)$ is said to be \begin{itemize}
	\item[---\,] \emph{geodesic} (or \emph{strictly intrinsic}) if, for any two $x,y \in X$, there exists a geodesic in $(X,d)$ from $x$ to $y$,
	\item[---\,] a \emph{length space} (or \emph{intrinsic}) if, for all $x,y \in X$ and $\epsilon \in \R_{>0}$, there exists a path in $(X,d)$ from $x$ to $y$ such that $\ell_{d}(\gamma) \leq d(x,y) + \epsilon$.
\end{itemize}

\begin{prop}[{\cite[\S2.4]{BuragoBuragoIvanov}}]\label{proposition:length.spaces} \begin{enumerate}
	\item\label{proposition:length.spaces.1} Every length space has approximate midpoints. 
	\item\label{proposition:length.spaces.2} Every complete metric space with approximate midpoints is a length space.
\end{enumerate} \end{prop}

\begin{proof} While~\ref{proposition:length.spaces.1} is~\cite[Lemma~2.4.10]{BuragoBuragoIvanov}, item~\ref{proposition:length.spaces.2} is~\cite[Theorem~2.4.16(2)]{BuragoBuragoIvanov}. \end{proof}

Of course, any length space is path-connected. Moreover, the following holds.

\begin{lem}\label{lemma:length.space.locally.path.connected} Every length space is locally path-connected. \end{lem}

\begin{proof} Let $(X,d)$ be a length space. In order to prove that its topology is locally path-connected, let $x \in X$ and $\epsilon \in \R_{>0}$. Consider the open ball $B \defeq \{ y \in X \mid d(x,y) < \epsilon \}$. Since $(X,d)$ is a length space, we see that \begin{displaymath}
	B \, = \, \bigcup \{ \im (\gamma) \mid \gamma \text{ path in } X, \, \gamma(0) = x, \, \ell_{d}(\gamma) < \epsilon \} .
\end{displaymath} In particular, $B$ is path-connected with respect to the induced topology. \end{proof}

A \emph{pointed metric space} is a triple $(X,d,x_{0})$ consisting of a metric space $(X,d)$ and a point $x_{0} \in X$. A pointed metric space $(X,d,x_{0})$ will be called \emph{star-shaped} if, for every $x \in X$, there exists a geodesic in $(X,d)$ from $x_{0}$ to $x$. A metric space $(X,d)$ is said to be \emph{star-shaped}~\cite{ThomWilson} if $(X,d,x_{0})$ is star-shaped for some~$x_{0} \in X$.

The following general construction involving geodesics will become relevant later, in the proof of Lemma~\ref{lemma:locally.algebraic}.

\begin{remark}\label{remark:geodesics.in.products} Let us consider a sequence of metric spaces $(X_{n},d_{n})_{n \in \N}$ such that $\sum\nolimits_{n \in \N} \diam(X_{n},d_{n}) < \infty$ and equip $X \defeq \prod_{n \in \N} X_{n}$ with the metric \begin{displaymath}
	d \colon \, X \times X \, \longrightarrow \, \R, \quad (x,y) \, \longmapsto \, \sum\nolimits_{n \in \N} d_{n}(x_{n},y_{n}) .
\end{displaymath} Let $x,y \in X$, define $s_{n} \defeq \sum\nolimits_{i=0}^{n-1} d_{i}(x_{i},y_{i})$ for every $n \in \N$, and note that $s_{0} = 0$. Now, for each $n \in \N$, let $\gamma_{n}$ be a geodesic from $x_{n}$ to $y_{n}$ in $(X_{n},d_{n})$. Then the map $\gamma \colon [0,d(x,y)] \to X$ defined by \begin{displaymath}
	\gamma(t)_{n} \, \defeq \, \begin{cases}
		\, x_{n} & \text{if } t \in [0,s_{n}) , \\
		\, \gamma_{n}(t-s_{n}) & \text{if } t \in [s_{n},s_{n+1}) , \\
		\, y_{n} & \text{otherwise, i.e., if } t \in [s_{n+1},d(x,y)]
	\end{cases}
\end{displaymath} for all $n \in \N$ and $t \in [0,d(x,y)]$ is a geodesic from $x$ to $y$ in $(X,d)$. \end{remark}

We briefly recollect some basic facts concerning metric ultraproducts, referring to~\cite[I.5]{BridsonHaefliger}, \cite[\S7]{RoeBook}, or~\cite[\S5]{BenYaacovBook} for further information. Let $(X_{i},d_{i})_{i \in I}$ be a family of metric spaces with $\sup_{i \in I} \diam(X_{i},d_{i}) < \infty$, and let $\mathcal{F}$ be an ultrafilter on the set $I$. Then the \emph{metric ultraproduct} $\prod\nolimits_{i \to \mathcal{F}} (X_{i},d_{i})$ is defined to be the metric quotient of $\prod\nolimits_{i \in I} X_{i}$ with respect to the pseudo-metric \begin{displaymath}
	\prod\nolimits_{i \in I} X_{i} \times \prod\nolimits_{i \in I} X_{i} \, \longrightarrow \, \R_{\geq 0}, \quad (x,y) \, \longmapsto \, \lim\nolimits_{i \to \mathcal{F}} d_{i}(x_{i},y_{i}) .
\end{displaymath} Recall that an ultrafilter $\mathcal{F}$ on a set $I$ is said to be \emph{countably incomplete} if there exists a countable subset $\mathcal{C} \subseteq \mathcal{F}$ such that $\bigcap \mathcal{C} = \emptyset$. As is well known and easy to see, any non-principal ultrafilter on a countable set is countably incomplete.

\begin{lem}\label{lemma:metric.ultraproduct} Let $\mathcal{F}$ be an ultrafilter on a set $I$, and let $(X_{i},d_{i})_{i \in I}$ be a family of metric spaces with $\sup_{i \in I} \diam(X_{i},d_{i}) < \infty$. Then the following hold. \begin{enumerate}
	\item\label{lemma:metric.ultraproduct.1} Consider a family of non-empty subsets $Y_{i} \subseteq X_{i}$ $(i \in I)$ such that \begin{displaymath}
			\qquad \lim\nolimits_{i \to \mathcal{F}}\sup\nolimits_{x \in X_{i}} d_{i}(x,Y_{i}) \, = \, 0 .
			\end{displaymath} If $\mathcal{F}$ is countably incomplete or $\{ i \in I \mid Y_{i} \textit{ closed in } (X_{i},d_{i})\} \in \mathcal{F}$, then \begin{displaymath}
			\qquad \iota \colon \, \prod\nolimits_{i \to \mathcal{F}} (Y_{i},d_{i}) \, \longrightarrow \, \prod\nolimits_{i \to \mathcal{F}} (X_{i},d_{i}), \quad [(y_{i})_{i \in I}] \, \longmapsto \, [(y_{i})_{i \in I}]
		\end{displaymath} is an isometric isomorphism.
	\item\label{lemma:metric.ultraproduct.2} If $\mathcal{F}$ is countably incomplete or $\{ i \in I \mid (X_{i},d_{i}) \textit{ complete}\} \in \mathcal{F}$, then the metric space $\prod\nolimits_{i \to \mathcal{F}} (X_{i},d_{i})$ is complete.
	\item\label{lemma:metric.ultraproduct.3} If $\mathcal{F}$ is countably incomplete and $\prod\nolimits_{i \to \mathcal{F}} (X_{i},d_{i})$ has approximate midpoints, then $\prod\nolimits_{i \to \mathcal{F}} (X_{i},d_{i})$ is geodesic.	
\end{enumerate} \end{lem}

\begin{proof} \ref{lemma:metric.ultraproduct.1} It is straightforward to verify that $\iota$ is a well-defined isometric embedding. To prove surjectivity, let $x \in \prod_{i \in I} X_{i}$. We proceed by case analysis.
	
\emph{Case 1}: $K \defeq \{ i \in I \mid Y_{i} \text{ closed in } (X_{i},d_{i}) \} \in \mathcal{F}$. Define $J \defeq \{ i \in I \mid x_{i} \in Y_{i} \}$. For each $i \in K \cap (I\setminus J)$, note that $d_{i}(x_{i},Y_{i}) > 0$ due to $Y_{i}$ being closed in $(X_{i},d_{i})$, and choose $y_{i} \in Y_{i}$ with $d_{i}(x_{i},y_{i}) \leq 2d_{i}(x_{i},Y_{i})$. For each $i \in (I\setminus K) \cup J$, set $y_{i} \defeq x_{i}$. Then $y \defeq (y_{i})_{i \in I} \in \prod_{i \in I} Y_{i}$ and \begin{displaymath}
	\lim\nolimits_{i \to \mathcal{F}} d_{i}(x_{i},y_{i}) \, \stackrel{K \in \mathcal{F}}{\leq} \, 2\lim\nolimits_{i \to \mathcal{F}} d_{i}(x_{i},Y_{i}) \, = \, 0 .
\end{displaymath} Thus, $[x] = \iota([y])$.	
	
\emph{Case 2}: $\mathcal{F}$ is countably incomplete. Then we find a descending chain $(I_{n})_{n \in \N}$ of members of $\mathcal{F}$ such that $\bigcap_{n \in \N} I_{n} = \emptyset$ and \begin{displaymath}
	\forall n \in \N \colon \quad I_{n} \, \subseteq \, \!\left\{ i \in I \left\vert \, d_{i}(x_{i},Y_{i}) < \tfrac{1}{n+1} \right\}\!\right. .
\end{displaymath} For each $i \in I_{0}$, since $(I_{n})_{n \in \N}$ is a descending chain with empty intersection, we may define $n(i) \defeq \max \{ n \in \N \mid i \in I_{n} \}$ and choose $y_{i} \in Y_{i}$ such that $d_{i}(x_{i},y_{i}) < \tfrac{1}{n(i)+1}$. For each $i \in I\setminus I_{0}$, pick any element $y_{i} \in Y_{i}$. Considering $y \defeq (y_{i})_{i \in I} \in \prod_{i \in I} Y_{i}$, we conclude that $[x] = \iota([y])$.
	
\ref{lemma:metric.ultraproduct.2} The case where $\{ i \in I \mid (X_{i},d_{i}) \text{ complete}\} \in \mathcal{F}$ is covered by~\cite[Proposition~5.3]{BenYaacovBook}. For the reader's convenience, we sketch the proof for the case where $\mathcal{F}$ is countably incomplete, which is a straightforward generalization of the argument for $I = \N$ given in~\cite[Lemma~5.53, p.~79]{BridsonHaefliger} and~\cite[Proposition~7.20, p.~108]{RoeBook}. Consider a sequence $x_{n} = (x_{n,i})_{i \in I} \in \prod_{i \in I} X_{i}$ $(n \in \N)$ such that $([x_{n}])_{n \in \N}$ is a Cauchy sequence in $\prod_{i \to \mathcal{F}} (X_{i},d_{i})$. Since $\mathcal{F}$ is countably incomplete, we find a descending chain $(I_{m})_{m \in \N}$ of members of $\mathcal{F}$ such that $\bigcap_{m \in \N} I_{m} = \emptyset$ and \begin{displaymath}
	\forall m \in \N \colon \quad I_{m} \, \subseteq \, \bigcap\nolimits_{n,n' \in \{ 0 ,\ldots,m \}} \!\left\{ i \in I \left\vert \, \vert d([x_{n}],[x_{n'}]) - d_{i}(x_{n,i},x_{n',i}) \vert \leq \tfrac{1}{m+1} \right\}\!\right. .
\end{displaymath} For each $i \in I_{0}$, since $(I_{m})_{m \in \N}$ is a descending chain with empty intersection, we may define $m(i) \defeq \max \{ m \in \N \mid i \in I_{m} \}$ as well as $y_{i} \defeq x_{m(i),i}$. Moreover, let $y_{i} \defeq x_{0,i}$ for all $i \in I\setminus I_{0}$. Considering $y \defeq (y_{i})_{i \in I}$, it is now easy to check that $([x_{n}])_{n \in \N}$ converges to $[y]$ in $\prod_{i \to \mathcal{F}} (X_{i},d_{i})$.

\ref{lemma:metric.ultraproduct.3} As $\mathcal{F}$ is countably incomplete, $\prod\nolimits_{i \to \mathcal{F}} (X_{i},d_{i})$ is complete by~\ref{lemma:metric.ultraproduct.2}. Since complete metric spaces with midpoints are necessarily geodesic according to~\cite[Theorem~2.4.16(1)]{BuragoBuragoIvanov}, thus it suffices to show that $\prod\nolimits_{i \to \mathcal{F}} (X_{i},d_{i})$ has midpoints. To this end, let $x,y \in \prod_{i \in I} X_{i}$ and put $t \defeq \tfrac{1}{2}d([x],[y])$. Since $\prod\nolimits_{i \to \mathcal{F}} (X_{i},d_{i})$ has approximate midpoints, we find a sequence $z_{n} \in \prod_{i \in I} X_{i}$ $(n \in \N)$ such that \begin{displaymath}
	\forall n \in \N \colon \quad \max\{\vert t - d([x],[z_{n}]) \vert, \vert t - d([z_{n}],[y]) \vert \} \, \leq \, \tfrac{1}{n+1} .
\end{displaymath} As $\mathcal{F}$ is countably incomplete, there exists a descending chain $(I_{n})_{n \in \N}$ of members of $\mathcal{F}$ such that $\bigcap_{n \in \N} I_{n} = \emptyset$ and \begin{displaymath}
	\forall n \in \N \colon \quad I_{n} \, \subseteq \, \!\left\{ i \in I \left\vert \, \max\{ \vert t - d_{i}(x_{i},z_{n,i}) \vert, \vert t - d_{i}(z_{n,i},y_{i}) \vert \} \leq \tfrac{1}{n+1} \right\}\!\right. .
\end{displaymath} For each $i \in I_{0}$, we define $n(i) \defeq \max \{ n \in \N \mid i \in I_{n} \}$ as well as $p_{i} \defeq z_{n(i),i}$. Moreover, let $p_{i} \defeq z_{0,i}$ for all $i \in I\setminus I_{0}$. Considering $p \defeq (p_{i})_{i \in I}$, we conclude that \begin{align*}
	d([x],[p]) \, &= \, \lim\nolimits_{i \to \mathcal{F}} d_{i}(x_{i},p_{i}) \, = \, t \, = \, \tfrac{1}{2}d([x],[y]), \\
	d([p],[y]) \, &= \, \lim\nolimits_{i \to \mathcal{F}} d_{i}(p_{i},y_{i}) \, = \, t \, = \, \tfrac{1}{2}d([x],[y]). \qedhere
\end{align*} \end{proof}

A \emph{metric group} (resp., \emph{pseudo-metric group}) is a pair $(G,d)$ consisting of a group~$G$ and a metric (resp., pseudo-metric) $d$ on $G$ such that \begin{displaymath}
	\forall g,x,y \in G \colon \qquad d(gx,gy) = d(x,y) = d(xg,yg) .
\end{displaymath} Note that, if $(G,d)$ is a pseudo-metric group, then $G$ constitutes a topological group with respect to the topology generated by $d$. Given any pseudo-metric group $(G,d)$ and $r \in \R_{>0}$, we let \begin{displaymath}
	\ball_{r}(G,d) \, \defeq \, \{ g \in G \mid d(1,g) \leq r \} .
\end{displaymath}

As is well known (see, e.g.,~\cite{ThomWilson,GMZ25}), if $\mathcal{F}$ is an ultrafilter on a set $I$ and $(G_{i},d_{i})_{i \in I}$ is a family of pseudo-metric groups with $\sup_{i \in I} \diam(G_{i},d_{i}) < \infty$, then the metric ultraproduct $\prod\nolimits_{i \to \mathcal{F}} (G_{i},d_{i})$ naturally constitutes a metric group, where \begin{displaymath}
	\prod\nolimits_{i \to \mathcal{F}} (G_{i},d_{i}) \, = \, \left(  \prod\nolimits_{i \in I} G_{i}\right)\!/N
\end{displaymath} for \begin{displaymath}
	\left. N \, \defeq \, \left\{ g \in \prod\nolimits_{i \in I} G_{i} \, \right\vert \lim\nolimits_{i \to \mathcal{F}} d_{i}(g_{i},1_{G_{i}}) = 0 \right\} \, \unlhd \, \prod\nolimits_{i \in I} G_{i} .
\end{displaymath}

We conclude this section with an observation about geodesics in metric groups.

\begin{definition} If $(G,d)$ is a metric group, then we define \begin{displaymath}
	\Delta(G,d) \, \defeq \, \{ g \in G \mid \exists \, \text{geodesic in $(X,d)$ from $1$ to $g$} \} .
\end{displaymath} \end{definition}

It is easy to see that a metric group $(G,d)$ is geodesic if and only if $\Delta(G,d) = G$. The following lemma examines a weaker hypothesis. 

\begin{lem}\label{lemma:approximate.midpoints} Let $(G,d)$ be a metric group with $G = \overline{\Delta(G,d)}$. Then $(G,d)$ has approximate midpoints. If $(G,d)$ is moreover complete, then $(G,d)$ is a length space. \end{lem}

\begin{proof} Let $g_{0},g_{1} \in G$ and $\epsilon \in \R_{>0}$. Since $G = \overline{\Delta(G,d)}$, there exists $h \in \Delta(G,d)$ such that $d(g_{0}g_{1}^{-1},h) \leq \tfrac{2\epsilon}{3}$. In particular, we find $s \in G$ with $d(1,s) = \tfrac{1}{2}d(1,h) = d(s,h)$. Considering $t \defeq sg_{1} \in G$, we see that \begin{align*}
	d(g_{0},t) \, &= \, d(g_{0}g_{1}^{-1},tg_{1}^{-1}) \, = \, d(g_{0}g_{1}^{-1},s) , \\
	d(g_{1},t) \, &= \, d(1,tg_{1}^{-1}) \, = \, d(1,s) \, = \, \tfrac{1}{2} d(h,1) ,
\end{align*} and therefore \begin{align*}
	\left\lvert d(g_{0},t) - \tfrac{1}{2}d(g_{0},g_{1}) \right\rvert \, &= \, \left\lvert d(g_{0}g_{1}^{-1},s) - \tfrac{1}{2}(g_{0}g_{1}^{-1},1) \right\rvert \\
		&\leq \, \left\lvert d(g_{0}g_{1}^{-1},s) - d(h,s) \right\rvert + \left\lvert \tfrac{1}{2} d(h,1) - \tfrac{1}{2} d(g_{0}g_{1}^{-1},1) \right\rvert \\
		&\leq \, \tfrac{3}{2} d(g_{0}g_{1}^{-1},h) \, \leq \, \epsilon , \\
		\left\lvert d(g_{1},t) - \tfrac{1}{2}d(g_{0},g_{1}) \right\rvert \, &= \, \left\lvert \tfrac{1}{2} d(h,1) - \tfrac{1}{2}(g_{0}g_{1}^{-1},1) \right\rvert \, \leq \, \tfrac{1}{2} d(g_{0}g_{1}^{-1},h) \, \leq \, \epsilon .
\end{align*} This shows that $(G,d)$ has approximate midpoints. The final assertion now follows due to Proposition~\ref{proposition:length.spaces}\ref{proposition:length.spaces.1}. \end{proof}

\section{Continuous rings}\label{section:continuous.rings}

The purpose of this second preliminary section is to put together some necessary background information regarding continuous rings and their rank functions. The reader is referred to von Neumann's original work~\cite{VonNeumannBook} and the monographs~\cite{MaedaBook,GoodearlBook} for a comprehensive account on the subject matter.

Let $R$ be a unital ring. Recall that the \emph{center} \begin{displaymath}
	\cent(R) \, \defeq \, \{ a \in R \mid \forall b \in R \colon \, ab = ba \}
\end{displaymath} is a commutative unital subring of $R$ and that the set \begin{displaymath}
	\GL(R) \, \defeq \, \{ a \in R \mid \exists b \in R \colon \, ab = ba = 1 \}
\end{displaymath} of \emph{units} of $R$ constitutes a group with respect to the multiplication inherited from~$R$. The set $\E(R) \defeq \{ e \in R \mid ee = e \}$ of \emph{idempotent} elements of $R$ admits a partial order defined by \begin{displaymath}
	e\leq f \quad :\Longleftrightarrow \quad ef=fe=e \qquad (e,f\in\E(R)).
\end{displaymath} We call two elements $e,f\in\E(R)$ \emph{orthogonal} and write $e\perp f$ if $ef=fe=0$. One readily checks that, for any $e,f \in \E(R)$, \begin{displaymath}
	e \perp f \quad \Longleftrightarrow \quad e \leq 1-f \quad \Longleftrightarrow \quad f \leq 1-e .
\end{displaymath}

A unital ring $R$ is said to be \emph{(von Neumann) regular}~\cite[II.II, Definition~2.2, p.~70]{VonNeumannBook} if, for every $a \in R$, there exists $b \in R$ such that $aba=a$. 

\begin{remark}[{\cite[II.II, Theorem~2.2, p.~70]{VonNeumannBook}}]\label{remark:regular} A unital ring $R$ is regular if and only if, for every $a \in R$, there exists $e \in \E(R)$ such that $aR=eR$. \end{remark}

A ring $R$ is called \emph{(directly) irreducible} if it is non-zero and not isomorphic to a direct product of two non-zero rings.

\begin{thm}[{\cite[II.II, Theorem~2.7, p.~75]{VonNeumannBook}}]\label{theorem:central.field} A regular ring $R$ is irreducible if and only if $\cent(R)$ is a field. \end{thm}

As established by von Neumann~\cite[II.II, Theorem~2.4, p.~72]{VonNeumannBook}, if $R$ is a regular ring, then the set $\lat(R) \defeq \{ aR \mid a \in R\}$, partially ordered by inclusion, constitutes a complemented, modular lattice. A \emph{continuous ring} is a regular ring $R$ such that $\lat(R)$ is a \emph{continuous geometry}~\cite[Chapter~13, p.~160--161]{GoodearlBook}, i.e., in addition to the aforementioned properties, the lattice $\lat(R)$ is complete and satisfies \begin{displaymath}
	I \wedge \bigvee C \, = \, \bigvee \{ I \wedge J \mid J \in C\}, \qquad \, I \vee \bigwedge C \, = \, \bigwedge \{ I \vee J \mid J \in C\} 
\end{displaymath} for every chain $C \subseteq \lat(R)$ and every $I \in \lat(R)$.

For our purposes, von Neumann's characterization of irreducible, continuous rings by means of the existence of a rank function (Theorem~\ref{theorem:unique.rank.function}) will be fundamental. In order to recollect some relevant terminology from~\cite[Chapter~16]{GoodearlBook}, let $R$ be a regular ring. A \emph{pseudo-rank function} on $R$ is a map $\rho \colon R \to [0,1]$ such that \begin{enumerate}
	\item[---\,] $\rho(1) = 1$,
	\item[---\,] $\rho(ab) \leq \min\{\rho(a),\rho(b)\}$ for all $a,b \in R$, and
	\item[---\,] $\rho(e+f) = \rho(e) + \rho(f)$ for any two orthogonal $e,f \in \E(R)$.
\end{enumerate} The last requirement implies that $\rho(0) = 0$ for any pseudo-rank function $\rho$ on $R$. A \emph{rank function} on $R$ is a pseudo-rank function $\rho$ on $R$ such that $\rho^{-1}(\{ 0 \}) = \{ 0 \}$.

\begin{lem}[\cite{VonNeumannBook}]\label{lemma:pseudo.rank.function} Let $\rho$ be a pseudo-rank function on a regular ring $R$. \begin{enumerate}
	\item\label{lemma:pseudo.rank.function.1} If $e,f \in \E(R)$ and $e \leq f$, then $\rho(f-e) = \rho(f)-\rho(e)$.
	\item\label{lemma:pseudo.rank.function.2} $\rho(a+b) \leq \rho(a) + \rho(b)$ for all $a,b \in R$.
\end{enumerate} \end{lem}

\begin{proof} \ref{lemma:pseudo.rank.function.1} For instance, see~\cite[Lemma~7.6(2)]{SchneiderGAFA}.

\ref{lemma:pseudo.rank.function.2} See~\cite[II.XVIII, Corollary~$(\overline{f})$ on p.~231]{VonNeumannBook}, \cite[VI.5, Hilfssatz~5.1(3°), p.~153]{MaedaBook}, or~\cite[Proposition~16.1(d)]{GoodearlBook}. \end{proof}

\begin{remark}[\cite{VonNeumannBook}]\label{remark:properties.pseudo.rank.function} Let $\rho$ be a pseudo-rank function on a regular ring $R$. \begin{enumerate}
	\item\label{remark:rank.ring.topology} It follows from Lemma~\ref{lemma:pseudo.rank.function}\ref{lemma:pseudo.rank.function.2} and the defining properties of a pseudo-rank function that \begin{align*}
			\qquad d_{\rho} \colon \, R \times R \, \longrightarrow \, [0,1], \quad (a,b) \, \longmapsto \, \rho(a-b)
		\end{align*} is a pseudo-metric on $R$, which is a metric if and only if $\rho$ is a rank function. Furnished with the \emph{$\rho$-topology}, i.e., the topology generated by $d_{\rho}$, the ring $R$ constitutes a topological ring (see~\cite[Remark~7.8]{SchneiderGAFA}). In particular, if $\rho$ is a rank function (equivalently, the $\rho$-topology is Hausdorff), then $\cent (R)$ is closed in $R$ with respect to the $\rho$-topology.
	\item\label{remark:unit.group} $(\GL(R),{d_{\rho}})$ is a pseudo-metric group (cf.~\cite[Remark~7.8]{SchneiderGAFA}). In particular, equipped with the relative $\rho$-topology, $\GL(R)$ is a topological group. Suppose now that $\rho$ is a rank function. Then \begin{displaymath}
		\qquad \GL(R) \, = \, \rho^{-1}(\{ 1\})
	\end{displaymath} is closed in $(R,{d_{\rho}})$ (see, e.g.,~\cite[Lemma~7.13(3)+(4)]{SchneiderGAFA}). Thus, if the metric space $(R,{d_{\rho}})$ is complete, then so is $(\GL(R),{d_{\rho}})$.
\end{enumerate} \end{remark}

\begin{thm}[\cite{VonNeumannBook}]\label{theorem:unique.rank.function} \begin{enumerate}
	\item\label{theorem:unique.rank.function.1} If $R$ is an irreducible, continuous ring, then $R$ admits a unique rank function $\rk_{R} \colon R \to [0,1]$, and the metric $d_{R} \defeq d_{\rk_{R}}$ is complete.
	\item\label{theorem:unique.rank.function.2} If $R$ is a regular ring admitting a rank function $\rho$ such that $(R,d_{\rho})$ is complete, then $R$ is a continuous ring.
\end{enumerate} \end{thm}

\begin{proof} \ref{theorem:unique.rank.function.1} Existence and uniqueness follow from~\cite[II.XVII, Theorem~17.1, p.~224]{VonNeumannBook} and~\cite[II.XVII, Theorem~17.2, p.~226]{VonNeumannBook}, respectively, while completeness is established in~\cite[II.XVII, Theorem~17.4, p.~230]{VonNeumannBook}.

\ref{theorem:unique.rank.function.2} Implicit already in~\cite[II.XVIII, Proof of Theorem~18.1, p.~237]{VonNeumannBook},
this is stated and proved explicitly in~\cite[VI.5, Satz~5.3, p.~156]{MaedaBook}. \end{proof}

For later use, we subsequently collect some structural properties reflected by the rank function of an irreducible, continuous ring. Given a unital ring $R$ and $n \in \N_{>0}$, we let $\M_{n}(R)$ denote the unital ring of $n \times n$ matrices with entries in $R$.

\begin{remark}[\cite{VonNeumannBook}]\label{remark:discrete} \begin{enumerate}
	\item\label{remark:discrete.1} Let $R$ be an irreducible, continuous ring. Due to work of von Neumann~\cite[I.VII, Theorem~7.3, p.~58]{VonNeumannBook} (see also~\cite[Remark~7.12]{SchneiderGAFA} or~\cite[Remark~3.4]{SchneiderIMRN}), the following are equivalent: \begin{enumerate}
			\item[---\,] $R$ is \emph{discrete}, i.e., the topology generated by~$d_{R}$ is discrete.
			\item[---\,] $\rk_{R}(R) = \left. \! \left\{ \tfrac{i}{n} \, \right\vert i \in \{ 0,\ldots,n \} \right\}$ for some $n \in \N_{>0}$.
			\item[---\,] $\rk_{R}(R) \ne [0,1]$.
		\end{enumerate}
	\item\label{remark:discrete.2} A ring $R$ is a discrete irreducible, continuous ring if and only if $R \cong \M_{n}(D)$ for some division ring $D$ and $n \in \N_{>0}$ (see~\cite[Remark~3.6]{SchneiderIMRN} for details). Moreover, if $D$ is a division ring and $n \in \N_{>0}$, then \begin{displaymath}
			\qquad \rk_{\M_{n}(D)}(a) \, = \, \tfrac{\rank(a)}{n}
		\end{displaymath} for all $a \in \M_{n}(D)$, as follows from Theorem~\ref{theorem:unique.rank.function}\ref{theorem:unique.rank.function.1} and the relevant properties of the natural rank map
	on $\M_{n}(D)$ (see~\cite[I.10.12, p.~359--360]{bourbakialg}).
\end{enumerate} \end{remark}

\begin{lem}[\cite{VonNeumannBook,GoodearlBook}]\label{lemma:rank.via.decompositions} Let $R$ be an irreducible, continuous ring. Moreover, let $x \in R$ and $n \in \N_{>0}$. Then \begin{displaymath}
	\rk_{R}(x) \geq \tfrac{1}{n} \ \ \Longleftrightarrow \ \ \exists a_{1},\ldots,a_{n},b_{1},\ldots,b_{n} \in R \colon \ 1 = \sum\nolimits_{i=1}^{n}a_{i}xb_{i}.
\end{displaymath} \end{lem}

\begin{proof} ($\Longleftarrow$) If $1 = \sum\nolimits_{i=1}^{n}a_{i}xb_{i}$ for some $a_{1},\ldots,a_{n},b_{1},\ldots,b_{n} \in R$, then \begin{displaymath}
	1 \, = \, \rk_{R}(1) \, = \, \rk_{R}\!\left( \sum\nolimits_{i=1}^{n}a_{i}xb_{i} \right)\! \, \stackrel{\ref{lemma:pseudo.rank.function}\ref{lemma:pseudo.rank.function.2}}{\leq} \, \sum\nolimits_{i=1}^{n} \rk_{R}(a_{i}xb_{i}) \, \leq \, n\rk_{R}(x) .
\end{displaymath}
	
($\Longrightarrow$) First, note that $\rk_{R}$ is the only pseudo-rank function on $R$. Indeed, since $R$ is simple by~\cite[VII.3, Hilfssatz~3.1, p.~166]{MaedaBook} (see also~\cite[Corollary~13.26, p.~170]{GoodearlBook}), any pseudo-rank function on $R$ is a rank function by~\cite[Proposition~16.7(a), p.~231]{GoodearlBook} and thus agrees with $\rk_{R}$ by Theorem~\ref{theorem:unique.rank.function}\ref{theorem:unique.rank.function.1}.
	
Now, let us suppose that $\rk_{R}(x) \geq \tfrac{1}{n}$, that is, $\rk_{R}(1) \leq n\rk_{R}(x)$. Then, according to~\cite[Theorem~18.28(a), p.~279]{GoodearlBook} and $\rk_{R}$ being the unique pseudo-rank function on~$R$, the right $R$-module $R_{R}$ embeds into the projective right $R$-module \begin{displaymath}
	M_{R} \, \defeq \, \underbrace{xR \oplus \ldots \oplus xR}_{\text{$n$-fold direct sum}}.
\end{displaymath} Since $R$ is regular, it follows by~\cite[Theorem~1.11, p.~6]{GoodearlBook} that $R_{R}$ is isomorphic to a direct summand of $M_{R}$. In particular, $R_{R}$ is a homomorphic image of $M_{R}$. Thanks to Remark~\ref{remark:regular}, there exists $e \in \E(R)$ such that $xR = eR$. As \begin{displaymath}
	\Hom(M_{R},R_{R}) \, \cong \, \Hom(xR_{R},R_{R}) \oplus \ldots \oplus \Hom(xR_{R},R_{R}) \, \cong \, Re \oplus \ldots \oplus Re ,
\end{displaymath} we find $a_{1},\ldots,a_{n} \in R$ with $1 \in \sum_{i=1}^{n}a_{i}xR$, which entails the desired conclusion. \end{proof}

A unital ring $R$ is said to be \emph{unit-regular}~\cite[Chapter~4, p.~37]{GoodearlBook} if, for every $a \in R$, there exists $u \in \GL(R)$ such that $aua = a$.

\begin{remark}[{\cite[Corollary~13.23, p.~170]{GoodearlBook}}]\label{remark:unit.regular} Every continuous ring is unit-regular. \end{remark}

\begin{lem}\label{lemma:invertible.approximation} Let $\rho$ be a pseudo-rank function on a unit-regular ring $R$. For every $a \in R$, there exists $b \in \GL(R)$ such that $\rho(b-a) = 1-\rho(a)$. \end{lem}

\begin{proof} Let $a \in R$. Since $R$ is unit-regular, there exists $u \in \GL(R)$ such that $aua = a$. Consider $b \defeq u^{-1} \in \GL(R)$. Then \begin{equation}\label{eq:equivalence}
	u(b-a) \, = \, ub-ua \, = \, 1-ua .
\end{equation} Since $uaua = ua$, i.e., $ua \in \E(R)$, we conclude that \begin{displaymath}
	\rho(b-a) \, \stackrel{u \in \GL(R)}{=} \, \rho(u(b-a)) \, \stackrel{\eqref{eq:equivalence}}{=} \, \rho(1-ua) \, \stackrel{\ref{lemma:pseudo.rank.function}\ref{lemma:pseudo.rank.function.1}}{=} \, 1-\rho(ua) \, \stackrel{u \in \GL(R)}{=} \, 1-\rho(a) . \qedhere
\end{displaymath} \end{proof}

Let us recall that the \emph{center} of a group $G$ is defined as \begin{displaymath}
	\cent (G) \, \defeq \, \{ h \in G \mid \forall g \in G \colon \, gh=hg \}
\end{displaymath} The center of the unit group of an irreducible, continuous ring is precisely the set of central units of the ring (Proposition~\ref{proposition:center}). This was established in~\cite[Lemma~8]{Ehrlich56} under the additional hypothesis that the ring's characteristic is different from $2$. Below we provide a proof for the general case: the main argument (Lemma~\ref{lemma:center}) follows an idea from~\cite[Lemma~2.2]{Halperin62}. We recall from~\cite[II.III, Definition~3.5, p.~97]{VonNeumannBook} the following piece of terminology concerning a unital ring $R$ and some $n \in \N_{>0}$: we say that $s \in \M_{n}(R)$ is a \emph{family of matrix units} for $R$ if $s_{11}+\ldots+s_{nn}=1$ and, for all $i,j,k,\ell\in\{1,\ldots,n\}$, \begin{displaymath}
	s_{ij}s_{k\ell} \, = \,
	\begin{cases}
		\, s_{i\ell} & \text{if } j=k, \\
		\, 0 & \text{otherwise}.
	\end{cases}
\end{displaymath}

\begin{lem}\label{lemma:center} Let $R$ be a unital ring admitting a family of matrix units in $\M_{m}(R)$ for some $m \in \N_{\geq 2}$. Then \begin{displaymath}
	\cent(R) \, = \, \{ a \in R \mid \forall b \in \GL(R) \colon \, ab = ba \} .
\end{displaymath} In particular, $\cent(\GL(R)) = \GL(R) \cap \cent(R)$. \end{lem}

\begin{proof} The proof is a variation of the argument in~\cite[Proof of Lemma~2.2]{Halperin62}.
	
The inclusion ($\subseteq$) is trivial. To prove ($\supseteq$), let $a \in R$ be such that $ab = ba$ for all $b \in \GL(R)$. Then \begin{equation}\label{partial.commute}
	\forall e \in \E(R) \ \forall x \in eR(1-e) \colon \qquad ax \, = \, xa .
\end{equation} Indeed, if $e \in \E(R)$ and $x \in eR(1-e)$, then $xx = x(1-e)ex = 0$ and therefore $(1+x)(1-x) = 1 = (1-x)(1+x)$, so that $1+x \in \GL(R)$ and hence our assumption about $a$ implies that \begin{displaymath}
	ax \, = \, a((1+x)-1) \, = \, a(1+x)-a \, = \, (1+x)a-a \, = \, ((1+x)-1)a \, = \, xa ,
\end{displaymath} as desired. Now, let $s \in \M_{m}(R)$ be a family of matrix units for $R$ with $m \in \N_{\geq 2}$. If $i,j \in \{ 1,\ldots,m \}$ are distinct, then the elements $s_{ii}, s_{jj} \in \E(R)$ are orthogonal, whence $s_{jj} \leq 1-s_{ii}$ and thus $s_{ii}Rs_{jj} \subseteq s_{ii}R(1-s_{ii})$, so that $ax=xa$ for all $x \in s_{ii}Rs_{jj}$ thanks to~\eqref{partial.commute}, and in turn \begin{displaymath}
	ax \, = \, as_{ii}x \, = \, as_{ij}s_{ji}x \, = \, s_{ij}as_{ji}x \, = \, s_{ij}s_{ji}xa \, = \, s_{ii}xa \, = \, xa
\end{displaymath} for all $x \in s_{ii}Rs_{ii}$. As $m \geq 2$, we conclude that \begin{displaymath}
	ax \, = \, \sum\nolimits_{i,j=1}^{m} as_{ii}xs_{jj} \, = \, \sum\nolimits_{i,j=1}^{m} s_{ii}xs_{jj}a \, = \, xa
\end{displaymath} for every $x \in R$, i.e., $a \in \cent(R)$. \end{proof}

\begin{prop}[cf.~{\cite[Lemma~8]{Ehrlich56}}]\label{proposition:center} Let $R$ be an irreducible, continuous ring. Then \begin{displaymath}
	\cent(R) \, = \, \{ a \in R \mid \forall b \in \GL(R) \colon \, ab = ba \} .
\end{displaymath} In particular, $\cent(\GL(R)) = \GL(R) \cap \cent(R) = \cent(R)\setminus \{ 0 \}$. \end{prop}

\begin{proof} If $R$ is a division ring, i.e., $\GL(R) = \cent(R)\setminus \{ 0 \}$, then the statement holds trivially. Henceforth, suppose that $R$ is not a division ring. Then, by Remark~\ref{remark:properties.pseudo.rank.function}\ref{remark:unit.group} and Remark~\ref{remark:discrete}\ref{remark:discrete.1}, there exists $m \in \N_{\geq 2}$ such that $\left. \! \left\{ \tfrac{i}{m} \, \right\vert i \in \{ 0,\ldots,m \} \right\} \subseteq \rk_{R}(R)$. Using the Hausdorff maximal principle, we find a maximal chain $E$ in $(\E(R),{\leq})$. According to~\cite[Corollary~7.19]{SchneiderGAFA}, we have $\rk_{R}(E) = \rk_{R}(R)$. Invoking Lemma~\ref{lemma:pseudo.rank.function}\ref{lemma:pseudo.rank.function.1}, we see that ${{\rk_{R}}\vert_{E}} \colon (E,{\leq}) \to (\rk_{R}(R),{\leq})$ is an isomorphism of linearly ordered sets. Since $\left. \! \left\{ \tfrac{i}{m} \, \right\vert i \in \{ 0,\ldots,m \} \right\} \subseteq \rk_{R}(R) = \rk_{R}(E)$, we may define \begin{displaymath}
	f_{i} \, \defeq \, ({\rk_{R}}\vert_{E})^{-1}\!\left(\tfrac{i}{m} \right) \qquad (i \in \{ 0,\ldots,m \}) .
\end{displaymath} Evidently, $0 = f_{0} \leq \ldots \leq f_{m} = 1$. Straightforward calculations show that \begin{displaymath}
	e_{i} \, \defeq \, f_{i} - f_{i-1} \, \in \, \E(R) \qquad (i \in \{ 1,\ldots,m \}) 
\end{displaymath} are pairwise orthogonal. Moreover, $\sum_{i=1}^{m} e_{i} = \sum_{i=1}^{m} f_{i} - \sum_{i=1}^{m} f_{i-1} = 1$. Also, for each $i \in \{1,\ldots,m\}$, \begin{displaymath}
	\rk_{R}(e_{i}) \, = \, \rk_{R}(f_{i} - f_{i-1}) \, \stackrel{\ref{lemma:pseudo.rank.function}\ref{lemma:pseudo.rank.function.1}}{=} \, \rk_{R}(f_{i}) - \rk_{R}(f_{i-1}) \, = \, \tfrac{i}{m} - \tfrac{i-1}{m} \, = \, \tfrac{1}{m} .
\end{displaymath} Hence, there exists a family of matrix units $s \in \M_{m}(R)$ for $R$ with $s_{ii} = e_{i}$ for each $i \in \{ 1,\ldots,m\}$ (see, e.g., \cite[Lemma~4.6]{BernardSchneider}). Since $m \geq 2$, we have \begin{displaymath}
	\cent(R) \, = \, \{ a \in R \mid \forall b \in \GL(R) \colon \, ab = ba \}
\end{displaymath} and $\cent(\GL(R)) = \GL(R) \cap \cent(R)$ by Lemma~\ref{lemma:center}. Finally, $\GL(R) \cap \cent(R) = \cent(R)\setminus \{ 0 \}$ thanks to Theorem~\ref{theorem:central.field}. \end{proof}

If $R$ is a unital ring and $e \in \E(R)$, then $eRe = \{eae \mid a \in R\}$ is a subring of $R$, with multiplicative identity $e$. The following facts concerning this standard construction for regular rings will be useful.

\begin{lem}[\cite{VonNeumannBook}]\label{lemma:corner.ring} Let $R$ be a regular ring and let $e \in \E(R)$. \begin{enumerate}
	\item\label{lemma:corner.ring.1} $eRe$ is a regular ring.
	\item\label{lemma:corner.ring.2} If $\rho$ is a rank function on $R$ and $e \ne 0$, then $\tfrac{1}{\rho(e)}\rho\vert_{eRe}$ constitutes a rank function on~$eRe$.
	\item\label{lemma:corner.ring.3} Suppose that $R$ is an irreducible, continuous ring and $e \ne 0$. Then $eRe$ is an irreducible, continuous ring with $\cent (eRe) = \cent(R)e \cong \cent(R)$ and \begin{displaymath}
			\qquad \rk_{eRe} \, = \, \tfrac{1}{\rk_{R}(e)}{\rk_{R}}\vert_{eRe} .
		\end{displaymath} Moreover, if $R$ is non-discrete, then so is $eRe$.
\end{enumerate} \end{lem}

\begin{proof} \ref{lemma:corner.ring.1} This is due to~\cite[II.II, Theorem~2.11, p.~77]{VonNeumannBook}.
	
\ref{lemma:corner.ring.2} This is straightforward (cf.~\cite[Lemma~16.2(a), p.~228]{GoodearlBook}).
	
\ref{lemma:corner.ring.3} The proof, based on references to~\cite{VonNeumannBook,Halperin62,GoodearlBook}, is given in~\cite[Remark~4.8]{BernardSchneider}. \end{proof}

We conclude this section by recalling some details concerning metric ultraproducts of irreducible, continuous rings.

\begin{prop}\label{proposition:ultraproduct} Let $(R_{i})_{i \in I}$ be a family of irreducible, continuous rings and let $\mathcal{F}$ be an ultrafilter on $I$. Then \begin{displaymath}
	\left. J \, \defeq \, \left\{ x \in \prod\nolimits_{i \in I} R_{i} \, \right\vert \lim\nolimits_{i \to \mathcal{F}} \rk_{R_{i}}(x_{i}) = 0 \right\}
\end{displaymath} is a two-sided ideal of the ring $\prod\nolimits_{i \in I} R_{i}$, and \begin{displaymath}
	\prod\nolimits_{i \to \mathcal{F}} R_{i} \, \defeq \, \!\left(  \prod\nolimits_{i \in I} R_{i}\right)\!/J
\end{displaymath} is an irreducible, continuous ring. Moreover, $\left(\prod\nolimits_{i \to \mathcal{F}} R_{i},d_{\prod\nolimits_{i \to \mathcal{F}} R_{i}}\right) = \prod_{i \to \mathcal{F}} (R_{i},d_{R_{i}})$. Finally, \begin{align*}
	\iota \colon \, \prod\nolimits_{i \to \mathcal{F}} (\GL(R_{i}),d_{R_{i}}) \, &\longrightarrow \, \left(\GL\!\left(\prod\nolimits_{i \to \mathcal{F}} R_{i}\right)\!,\,d_{\prod\nolimits_{i \to \mathcal{F}} R_{i}}\right)\!, \\
	[(x_{i})_{i \in I}] \, &\longmapsto \, [(x_{i})_{i \in I}]
\end{align*} is an isometric group isomorphism. \end{prop}

\begin{proof} Using the properties of a rank function, it is easy to verify that $J$ is indeed a two-sided ideal of $\prod\nolimits_{i \in I} R_{i}$. Since the class of regular rings is closed under direct products and quotients, $R \defeq \prod\nolimits_{i \to \mathcal{F}} R_{i}$ is a regular ring. Note that $R$ is non-zero, as $1_{R} \notin J$. Straightforward calculations reveal that $\rho \colon R \to [0,1], \, [(x_{i})_{i \in I}] \mapsto \lim\nolimits_{i \to \mathcal{F}} \rk_{R_{i}}(x_{i})$ constitutes a well-defined rank function on $R$ and that $(R,d_{\rho}) = \prod_{i \to \mathcal{F}} (R_{i},d_{R_{i}})$. In particular, $(R,d_{\rho})$ is complete by Theorem~\ref{theorem:unique.rank.function}\ref{theorem:unique.rank.function.1} and Lemma~\ref{lemma:metric.ultraproduct}\ref{lemma:metric.ultraproduct.2}, and hence $R$ is a continuous ring by Theorem~\ref{theorem:unique.rank.function}\ref{theorem:unique.rank.function.2}. Since $R$ is non-zero, irreducibility will follow from simplicity. We will prove the simplicity of $R$ by suitably adapting an argument from~\cite[Proof of Theorem~4.6, p.~65]{stolz}. To this end, consider any $x = [(x_{i})_{i \in I}] \in R\setminus \{ 0 \}$. Then $\tau \defeq \lim_{i\to \mathcal{F}}\rk_{R_{i}}(x_{i}) > 0$. Fix any $n \in \N_{>0}$ such that $1/n < \tau$. Then \begin{displaymath}
	I_{0} \, \defeq \, \!\left\{ i \in I \left\vert \, \rk_{R_{i}}(x_{i}) > \tfrac{1}{n} \right\} \!\right. \, \in \, \mathcal{F} .
\end{displaymath} For each $i \in I_{0}$, we invoke Lemma~\ref{lemma:rank.via.decompositions} to find $a_{1,i},\ldots,a_{n,i},b_{1,i},\ldots,b_{n,i} \in R_{i}$ such that $1_{R_{i}} = \sum\nolimits_{k=1}^{n} a_{k,i}x_{i}b_{k,i}$. For all $i \in I\setminus I_{0}$ and $k \in \{ 1,\ldots,n\}$, we let $a_{k,i} \defeq 0_{R_{i}}$ and $b_{k,i} \defeq 0_{R_{i}}$. Considering the elements \begin{displaymath}
	a_{k} \defeq [(a_{k,i})_{i \in I}] \in R, \qquad b_{k} \defeq [(b_{k,i})_{i \in I}] \in R \qquad (k \in \{ 1,\ldots,n \}),
\end{displaymath} we conclude that $1_{R} = \sum_{k=1}^{n} a_{k}xb_{k}$. It follows that every non-zero two-sided ideal of $R$ coincides with $R$. This means that $R$ is simple, thus irreducible. Hence, $\rho = \rk_{R}$ by Theorem~\ref{theorem:unique.rank.function}\ref{theorem:unique.rank.function.1}. Now, it is easy to see that $\iota$ is a well-defined isometric group embedding. To prove that $\iota$ is surjective, let $x \in \prod\nolimits_{i \in I} R_{i}$ with $[x] \in \GL(R)$. Then \begin{displaymath}
	1 \, \stackrel{\ref{remark:properties.pseudo.rank.function}\ref{remark:unit.group}}{=} \, \rk_{R}([x]) \, = \, \lim\nolimits_{i \to \mathcal{F}} \rk_{R_{i}}(x_{i}) 
\end{displaymath} and therefore $I'_{n} \defeq \!\left\{ i \in I \left\vert \, \rk_{R_{i}}(x_{i}) \geq 1 - \tfrac{1}{n} \right\}\!\right. \in \mathcal{F}$ for every $n \in \N_{>0}$. We observe that $(I'_{n})_{n \in \N_{>0}}$ constitutes a descending chain and consider \begin{displaymath}
	I'_{\infty} \, \defeq \, \bigcap\nolimits_{n \in \N_{>0}} I'_{n} \, = \, \{ i \in I \mid \rk_{R_{i}}(x_{i}) = 1 \} .
\end{displaymath} For each $i \in I'_{1}\setminus I'_{\infty}$, we let $n(i) \defeq \max \{ n \in \N_{>0} \mid i \in I'_{n} \}$ and find, using Remark~\ref{remark:unit.regular} and Lemma~\ref{lemma:invertible.approximation}, an element $y_{i} \in \GL(R_{i})$ such that $d_{R_{i}}(x_{i},y_{i}) < \tfrac{1}{n(i)}$. Moreover, let $y_{i} \defeq x_{i}$ for all $i \in I'_{\infty}$ and $y_{i} \defeq 1_{R_{i}}$ for all $i \in I\setminus I'_{1}$. Considering $y \defeq (y_{i})_{i \in I}$, we observe that $[y] \in \prod\nolimits_{i \to \mathcal{F}} (\GL(R_{i}),d_{R_{i}})$. Finally, as $I'_{n} \in \mathcal{F}$ for each $n \in \N_{>0}$, \begin{displaymath}
	d_{R}([x],\iota([y])) \, = \, \lim\nolimits_{i \to \mathcal{F}} d_{R_{i}}(x_{i},y_{i}) \, = \, 0 
\end{displaymath} and thus $[y] =\iota([x])$, as desired. \end{proof}

\section{Topological simplicity}\label{section:topological.simplicity}

The purpose of this section is the proof of Theorem~\ref{theorem:tbng} and Corollary~\ref{corollary:ultraproduct}. Our argument rests on an application of a result due to Rodgers and Saxl (Theorem~\ref{theorem:rodgers.saxl}), the statement of which involves a notion of \emph{index} (Definition~\ref{definition:index}) related to the rational canonical form of a square matrix (Theorem~\ref{theorem:rational.canonical.form}). 

To clarify some notation, let $K$ be a field. Consider the (univariate) polynomial ring $K[X]$ over $K$ along with the usual degree function $\deg \colon K[X]\setminus \{ 0 \} \to \N$. The \emph{companion matrix} of a monic polynomial $f = X^{n} + \sum_{i=0}^{n-1} a_{i}X^{i} \in K[X]$ with $n \defeq \deg(f) > 0$ is defined as \begin{displaymath}
	\comp(f) \, \defeq \, \begin{pmatrix}
		0 & \ldots & \ldots & 0 & -a_{0} \\
		1 & \ddots & & \vdots & \vdots \\
		0 & \ddots & \ddots & \vdots & \vdots \\
		\vdots & \ddots & \ddots & 0 & \vdots \\
		0 & \ldots	& 0 & 1 & -a_{n-1}	
	\end{pmatrix} \, \in \, \M_{n}(K) .
\end{displaymath}

\begin{thm}[{\cite[V.21.4, Theorem~4.4, p.~417]{linalg}}]\label{theorem:rational.canonical.form} Let $n \in \N_{>0}$, let $K$ be a field, and let $a \in \M_{n}(K)$. Then there exists a unique tuple \begin{displaymath}
	(f_{1},\ldots,f_{r}) \, \in \, \bigcup\nolimits_{s \in \N_{>0}} K[X]^{s} ,
\end{displaymath} called the \emph{rational canonical form} of $a$, such that \begin{itemize}
	\item[---\,] $\deg(f_{1}) > 0$,
	\item[---\,] $f_{1},\ldots,f_{r}$ are monic,
	\item[---\,] $f_{i} \mid f_{i+1}$ for each $i \in \{ 1,\ldots, r-1\}$, and
	\item[---\,] there exists $b \in \GL_{n}(K)$ with \begin{displaymath}
					\qquad bab^{-1} \, = \, \begin{pmatrix} 
											\comp(f_{1}) & & 0 \\
											& \ddots & \\
											0 & & \comp(f_{r})
										\end{pmatrix} .
				\end{displaymath}
\end{itemize} \end{thm}

\begin{definition}[{\cite[pp.~4625--4626]{RodgersSaxl}}; see also~{\cite[Section~2]{KnueppelNielsen}}]\label{definition:index} Let $K$ be a field and let $n \in \N_{>0}$. If $(f_{1},\ldots,f_{r})$ is the rational canonical form of $a \in \M_{n}(K)$, then\footnote{In~\cite{RodgersSaxl,KnueppelNielsen}, this definition is given only for invertible square matrices.} \begin{displaymath}
	\ind(a) \, \defeq \, n-r
\end{displaymath} is called the \emph{index} of $a$. \end{definition}

\begin{thm}[{\cite[Theorem~2.8]{RodgersSaxl}}]\label{theorem:rodgers.saxl} Let $K$ be a field, $n \in \N_{>2}$, $a_{1},\ldots,a_{m} \in \SL_{n}(K)$. If $\sum\nolimits_{i=1}^{m} \ind(a_{i}) > 6(n-1)$, then $\SL_{n}(K) = \Cl_{\SL_{n}(K)}(a_{1})\cdots \Cl_{\SL_{n}(K)}(a_{m})$.  \end{thm}

The following simple lemma is of central importance to the proof of Theorem~\ref{theorem:tbng}, as it relates the index of a square matrix to the corresponding rank distance from the center and thus facilitates the application of Theorem~\ref{theorem:rodgers.saxl}.

\begin{lem}\label{lemma:index} Let $n \in \N_{>0}$, let $K$ be a field, let $a \in \M_{n}(K)$. Then \begin{displaymath}
	\inf\nolimits_{c \in K} \rank (a-c1) \, \leq \, 2 \ind (a) .
\end{displaymath} Equivalently, \begin{displaymath}
	d_{\M_{n}(K)}(a,\cent(\M_{n}(K))) \, \leq \, \tfrac{2}{n} \ind (a) .
\end{displaymath} \end{lem}

\begin{proof} Let $(f_{1},\ldots,f_{r})$ be the rational canonical form of $a$. Define \begin{displaymath}
	t \, \defeq \, \max (\{ 0 \} \cup \{ i \in \{1,\ldots,r\} \mid \deg(f_{i}) = 1 \} ) .
\end{displaymath} Since $f_{1},\ldots,f_{t}$ are monic and $f_{1} \mid f_{2} \mid \ldots \mid f_{t}$, we see that $f_{1} = \ldots = f_{t} = X - c$ for some $c \in K$. Furthermore, \begin{displaymath}
	n \, = \, \sum\nolimits_{i=1}^{r} \deg(f_{i}) \, = \, t + \sum\nolimits_{i=t+1}^{r} \deg(f_{i}) \, \geq \, t + 2(r-t) \, = \, 2r-t
\end{displaymath} and hence \begin{displaymath}
	n-t \, \leq \, 2(n-r) \, = \, 2\ind(a) .
\end{displaymath} Finally, we find $b \in \GL_{n}(K)$ such that \begin{displaymath}
	bab^{-1} \, = \, \begin{pmatrix} 
	\comp(f_{1}) & & 0 \\
		& \ddots & \\
		0 & & \comp(f_{r})
	\end{pmatrix} .
\end{displaymath} As $\comp(f_{i}) = \comp(X-c) = (c) \in \M_{1}(K)$ for each $i \in \{ 1,\ldots,t\}$, we conclude that \begin{align*}
	\rank(a-c1) \, &= \, \rank\!\left(b(a-c1)b^{-1}\right) \, = \, \rank\!\left(bab^{-1}-c1\right) \\
	& = \, \rank \left( \begin{pmatrix} \comp(f_{1}) & & 0 \\
		& \ddots & \\
		0 & & \comp(f_{r}) \end{pmatrix} - \begin{pmatrix} c & & 0 \\
		& \ddots & \\
		0 & & c \end{pmatrix} \right) \\
	&= \, \rank \begin{pmatrix}
			0 & & & & & 0 \\
			& \ddots & & & & \\
			& & 0 & & & \\
			& & & \comp(f_{t+1})-c & &  \\
			& & & & \ddots & \\
			0 & & & & & \comp(f_{r})-c \end{pmatrix} \\
	&\leq \, n-t \, \leq \, 2 \ind(a) .
\end{align*} Due to Remark~\ref{remark:discrete}\ref{remark:discrete.2} and the fact that $\cent(\M_{n}(K)) = \{ c1 \mid c \in K \}$, the final statement is an equivalent reformulation. \end{proof}

We now arrive at a reinforced extension to arbitrary fields of~\cite[Theorem~3.1]{CarderiThom} (see also~\cite[Lemma~9.1]{GMZ25}). The latter constitutes a consequence of work of Liebeck and Shalev~\cite[Lemma~5.4]{LiebeckShalev}.

\begin{cor}\label{corollary:index} Let $n \in \N_{>2}$, let $K$ be a field, and let $a_{1},\ldots,a_{m} \in \SL_{n}(K)$. If $\sum_{i=1}^{m} d_{\M_{n}(K)}(a_{i},\cent(\M_{n}(K))) \geq 12$, then $\SL_{n}(K) = \Cl_{\SL_{n}(K)}(a_{1})\cdots \Cl_{\SL_{n}(K)}(a_{m})$. \end{cor}

\begin{proof} Since \begin{displaymath}
	\sum\nolimits_{i=1}^{m} \ind(a_{i}) \, \stackrel{\ref{lemma:index}}{\geq} \, \tfrac{n}{2} \sum\nolimits_{i=1}^{m} d_{\M_{n}(K)}(a_{i},\cent(\M_{n}(K))) \, \geq \, 6n \, > \, 6(n-1) ,
\end{displaymath} the claim follows by Theorem~\ref{theorem:rodgers.saxl}. \end{proof}

Another important ingredient in the proof of Theorem~\ref{theorem:tbng} is the following approximation theorem from~\cite{BernardSchneider}.

\begin{prop}[{\cite[Proposition~9.6]{BernardSchneider}}]\label{proposition:simply.special.dense} Let $R$ be a non-discrete irreducible, continuous ring, $K \defeq \cent(R)$, $g \in \GL(R)$ and $\epsilon \in \R_{>0}$. Then there exists $m \in \N_{>0}$ such that, for every $n \in \N_{>0}$ with $m \vert n$, there exist a unital $K$-algebra embedding $\phi \colon \M_{n}(K) \to R$ and an element $a \in \SL_{n}(K)$ such that $d_{R}(g,\phi(a)) < \epsilon$. \end{prop}

In order to combine multiple embeddings of a matrix algebra, we will need the following observation, a variant of the Skolem--Noether theorem (see, e.g.,~\cite[I.3, Theorem~3.14, p.~93]{FarbDennis}). Given a group $G$ and $g \in G$, we let $\gamma_{g} \colon G \to G, \, x \mapsto gxg^{-1}$.

\begin{lem}\label{lemma:conjugation} Let $R$ be an irreducible, continuous ring, let $K \defeq \cent (R)$, let $n \in \N_{>0}$, and let $\phi, \psi \colon \M_{n}(K) \to R$ be two unital $K$-algebra embeddings. Then there exists $g \in \GL(R)$ such that $\psi = {\gamma_{g}} \circ \phi$. \end{lem}

\begin{proof} By a standard fact (see, e.g.,~\cite[Remark~3.10]{BernardSchneider}), there exist two families of matrix units $s,t \in \M_{n}(R)$ for $R$ such that, for all $a \in \M_{n}(K)$, \begin{displaymath}
	\phi(a) \, = \, \sum\nolimits_{i,j=1}^{n} a_{ij}s_{ij}, \qquad \psi(a) \, = \, \sum\nolimits_{i,j=1}^{n} a_{ij}t_{ij}.
\end{displaymath} Note that $\rk_{R}(s_{11}) = \tfrac{1}{n} = \rk_{R}(t_{11})$ by~\cite[Remark~4.1(B)]{BernardSchneider}. Thus, due to~\cite[Lemma~9]{Ehrlich56}, there exists $u \in \GL(R)$ such that $t_{11} = us_{11}u^{-1}$, i.e., $s_{11} = u^{-1}t_{11}u$. For each $i \in \{ 1,\ldots,n \}$, define \begin{displaymath}
	v_{i} \, \defeq \, t_{i1}us_{1i} \, \in \, t_{ii}Rs_{ii}, \qquad w_{i} \, \defeq \, s_{i1}u^{-1}t_{1i} \, \in \, s_{ii}Rt_{ii} .
\end{displaymath} For any two distinct $i,j \in \{ 1,\ldots,n \}$, since $s_{1i}s_{j1} = 0$ and $t_{1j}t_{i1} = 0$, we see that \begin{equation}\label{eq:orthogonal}
	v_{i}w_{j} \, = \, t_{i1}us_{1i}s_{j1}u^{-1}t_{1j} \, = \, 0 \, = \, s_{j1}u^{-1}t_{1j}t_{i1}us_{1i} \, = \, w_{j}v_{i} .
\end{equation} Now, consider $g \defeq v_{1} + \ldots + v_{n}$ and $h \defeq w_{1} + \ldots + w_{n}$. Then \begin{align*}
	gh \, &= \, (v_{1} + \ldots + v_{n})(w_{1} + \ldots + w_{n}) \, \stackrel{\eqref{eq:orthogonal}}{=} \, v_{1}w_{1} + \ldots + v_{n}w_{n} \\
	& = \, t_{11}us_{11}s_{11}u^{-1}t_{11} + \ldots + t_{n1}us_{1n}s_{n1}u^{-1}t_{1n} \\
	& = \, t_{11}us_{11}u^{-1}t_{11} + \ldots + t_{n1}us_{11}u^{-1}t_{1n} \, = \, t_{11}t_{11}t_{11} + \ldots + t_{n1}t_{11}t_{1n} \\
	& = \, t_{11} + \ldots + t_{nn} \, = \, 1, \\
	hg \, &= \, (w_{1} + \ldots + w_{n})(v_{1} + \ldots + v_{n}) \, \stackrel{\eqref{eq:orthogonal}}{=} \, w_{1}v_{1} + \ldots + w_{n}v_{n} \\
	& = \, s_{11}u^{-1}t_{11}t_{11}us_{11} + \ldots + s_{n1}u^{-1}t_{1n}t_{n1}us_{1n} \\
	& = \, s_{11}u^{-1}t_{11}us_{11} + \ldots + s_{n1}u^{-1}t_{11}us_{1n} \, = \, s_{11}s_{11}s_{11} + \ldots + s_{n1}s_{11}s_{1n} \\
	& = \, s_{11} + \ldots + s_{nn} \, = \, 1.
\end{align*} This shows that $g \in \GL(R)$ with $g^{-1} = h$. Finally, for each pair $(i,j) \in \{ 1,\ldots,n\}^{2}$, \begin{align*}
	gs_{ij}g^{-1} \, &= \, gs_{ij}h \, = \, (v_{1} + \ldots + v_{n})s_{ij}(w_{1} + \ldots + w_{n}) \, = \, v_{i}s_{ij}w_{j} \\
	& = \, t_{i1}us_{1i}s_{ij}s_{j1}u^{-1}t_{1j} \, = \, t_{i1}us_{11}u^{-1}t_{1j} \, = \, t_{i1}t_{11}t_{1j} \, = \, t_{ij} .
\end{align*} This entails that \begin{displaymath}
	\psi(a) \, = \, \sum\nolimits_{i,j=1}^{n} a_{ij}t_{ij} \, = \, \sum\nolimits_{i,j=1}^{n} a_{ij}gs_{ij}g^{-1} \, = \, g\!\left( \sum\nolimits_{i,j=1}^{n} a_{ij}s_{ij} \right)\!g^{-1} \, = \, \gamma_{g}(\phi(a))
\end{displaymath} for all $a \in \M_{n}(K)$, i.e., $\psi = {\gamma_{g}} \circ \phi$. \end{proof}

We are now ready to prove this section's main result.

\begin{thm}\label{theorem:tbng} Let $R$ be a non-discrete irreducible, continuous ring and let $K \defeq \cent (R)$. If $m \in \N$ and $g_{1},\ldots,g_{m} \in \GL(R)$ are such that $\sum\nolimits_{i=1}^{m} d_{R}(g_{i},K) > 12$, then \begin{displaymath}
	\GL(R) \, = \, \overline{\Cl_{\GL(R)}(g_{1})\cdots \Cl_{\GL(R)}(g_{m})} .
\end{displaymath} In particular, if $g \in \GL(R)\setminus K$ and $m \in \N$ are such that\footnote{Since $K$ is closed in $(R,d_{R})$ by Remark~\ref{remark:properties.pseudo.rank.function}\ref{remark:rank.ring.topology}, we have $d_{R}(g,K) > 0$ for every $g \in R\setminus K$.} $m > \tfrac{12}{d_{R}(g,K)}$, then \begin{displaymath}
	\GL(R) \, = \, \overline{\Cl_{\GL(R)}(g)^{m}} .
\end{displaymath} \end{thm}

\begin{proof} Let $m \in \N$ and $g_{1},\ldots,g_{m} \in \GL(R)$ such that $\tau \defeq \left(\sum\nolimits_{i=1}^{m} d_{R}(g_{i},K)\right) - 12 > 0$. Consider any $\tilde{g} \in \GL(R)$. We will prove that $\tilde{g} \in \overline{\Cl_{\GL(R)}(g_{1})\cdots \Cl_{\GL(R)}(g_{m})}$. To this end, let $\epsilon \in \R_{>0}$. By Proposition~\ref{proposition:simply.special.dense}, there exist $n \in \N_{>2}$, unital $K$-algebra embeddings $\phi_{1},\ldots,\phi_{m}, \tilde{\phi} \colon \M_{n}(K) \to R$ and $a_{1},\ldots,a_{m},\tilde{a} \in \SL_{n}(K)$ such that \begin{displaymath}
	\forall i \in \{ 1,\ldots,m\} \colon \qquad d_{R}(g_{i},\phi_{i}(a_{i})) \, \leq \, \tfrac{\min \{\epsilon,\tau \}}{2m}
\end{displaymath} and $d_{R}(\tilde{g},\tilde{\phi}(\tilde{a})) \leq \tfrac{\epsilon}{2}$. Thanks to Lemma~\ref{lemma:conjugation}, there exist $h_{1},\ldots,h_{m} \in \GL(R)$ such that $\tilde{\phi} = {\gamma_{h_{i}}} \circ \phi_{i}$ for each $i \in \{ 1,\ldots,m \}$. Since $\M_{n}(K)$ is an irreducible, continuous ring by Remark~\ref{remark:discrete}\ref{remark:discrete.2}, we know from~\cite[Remark~4.5(B)]{BernardSchneider} that $d_{\M_{n}(K)}(a,b) = d_{R}(\phi_{i}(a),\phi_{i}(b))$ for all $a,b \in \M_{n}(K)$ and $i \in \{ 1,\ldots,m\}$. In particular, \begin{align*}
	\sum\nolimits_{i=1}^{m} d_{\M_{n}(K)}(a_{i},\cent(\M_{n}(K))) \, &= \, \sum\nolimits_{i=1}^{m} d_{R}(\phi_{i}(a_{i}),\phi_{i}(\cent(\M_{n}(K)))) \\
	& \geq \, \left(\sum\nolimits_{i=1}^{m} d_{R}(g_{i},K) \right) - m \cdot \tfrac{\tau}{2m} \, = \, 12+\tau -\tfrac{\tau}{2} \, > \, 12 .
\end{align*} Consequently, due to Corollary~\ref{corollary:index}, we find $b_{1},\ldots,b_{m} \in \SL_{n}(K)$ such that \begin{displaymath}
	\tilde{a} \, = \, b_{1}a_{1}b_{1}^{-1}\cdots b_{m}a_{m}b_{m}^{-1} .
\end{displaymath} For each $i \in \{ 1,\ldots,m \}$, consider \begin{align*}
	\tilde{g}_{i} \, &\defeq \, h_{i}\phi_{i}(b_{i})g_{i}\phi_{i}(b_{i})^{-1}h_{i}^{-1} \, \in \, \Cl_{\GL(R)}(g_{i}), \\
	\dtilde{g}_{i} \, &\defeq \, h_{i}\phi_{i}(b_{i})\phi_{i}(a_{i})\phi_{i}(b_{i})^{-1}h_{i}^{-1} \, \in \, \GL(R).
\end{align*} We observe that \begin{align*}
	\tilde{\phi}(\tilde{a}) \, &= \, \tilde{\phi}\!\left(b_{1}a_{1}b_{1}^{-1}\cdots b_{m}a_{m}b_{m}^{-1} \right)\! \, = \, \tilde{\phi}\!\left(b_{1}a_{1}b_{1}^{-1}\right) \cdots \tilde{\phi}\!\left(b_{m}a_{m}b_{m}^{-1} \right)\\
	& = \,h_{1}\phi_{1}\!\left(b_{1}a_{1}b_{1}^{-1}\right)\!h_{1}^{-1} \cdots h_{m}\phi_{m}\!\left(b_{m}a_{m}b_{m}^{-1} \right)\!h_{m}^{-1} \, = \, \dtilde{g}_{1}\cdots \dtilde{g}_{m}
\end{align*} and hence \begin{align*}
	d_{R}(\tilde{g},\tilde{g}_{1}\cdots \tilde{g}_{m}) \, &\leq \, d_{R}(\tilde{g},\tilde{\phi}(\tilde{a})) + d_{R}(\tilde{\phi}(\tilde{a}),\tilde{g}_{1}\cdots \tilde{g}_{m}) \, \leq \, \tfrac{\epsilon}{2} + d_{R}(\dtilde{g}_{1}\cdots \dtilde{g}_{m},\tilde{g}_{1}\cdots \tilde{g}_{m}) \\
	&\leq \, \tfrac{\epsilon}{2} + \sum\nolimits_{i=1}^{m} d_{R}(\tilde{g}_{1}\cdots \tilde{g}_{i-1}\dtilde{g}_{i}\dtilde{g}_{i+1}\cdots \dtilde{g}_{m},\tilde{g}_{1}\cdots \tilde{g}_{i-1}\tilde{g}_{i}\dtilde{g}_{i+1}\cdots \dtilde{g}_{m}) \\
	& \stackrel{\ref{remark:properties.pseudo.rank.function}\ref{remark:unit.group}}{=} \, \tfrac{\epsilon}{2} + \sum\nolimits_{i=1}^{m} d_{R}(\dtilde{g}_{i},\tilde{g}_{i}) \, \stackrel{\ref{remark:properties.pseudo.rank.function}\ref{remark:unit.group}}{=} \, \tfrac{\epsilon}{2} +\sum\nolimits_{i=1}^{m} d_{R}(\phi_{i}(a_{i}),g_{i}) \\
	& \leq \, \tfrac{\epsilon}{2} + m \cdot \tfrac{\epsilon}{2m} \, = \, \epsilon .
\end{align*} This shows that $\tilde{g} \in \overline{\Cl_{\GL(R)}(g_{1})\cdots \Cl_{\GL(R)}(g_{m})}$, as intended. The final assertion readily follows. \end{proof}

To spell out further consequences of Theorem~\ref{theorem:tbng}, we recall the following.

\begin{definition} A group $G$ has \emph{bounded normal generation}~\cite[Definition~4.7(ii)]{Dowerk} (see also~\cite[Definition~2.1]{DowerkLeMaitre}) if \begin{displaymath}
	\forall g \in G\setminus \{ e\} \ \exists n \in \N \colon \qquad G \, = \, \!\left(\Cl_{G}(g)\cup \Cl_{G}\!\left(g^{-1}\right)\right)^{n} .
\end{displaymath} A topological group $G$ has \emph{topological bounded normal generation}~\cite[Definition~4.7(iv)]{Dowerk} (see also~\cite[Definition~2.2]{DowerkLeMaitre}) if \begin{displaymath}
	\forall g \in G\setminus \{ e\} \ \exists n \in \N \colon \qquad G \, = \, \overline{\left(\Cl_{G}(g)\cup \Cl_{G}\!\left(g^{-1}\right)\right)^{n}} .
\end{displaymath} \end{definition}

\begin{remark}[{\cite[Proposition~4.8(i)]{Dowerk}}]\label{remark:topological.simplicity} Any group with bounded normal generation is simple. Any topological group with topological bounded normal generation is \emph{topologically simple}, i.e., it does not admit any non-trivial closed normal subgroup. \end{remark}

\begin{cor}\label{corollary:pgl} Let $R$ be a non-discrete irreducible, continuous ring. Then $\PGL(R)$ has topological bounded normal generation. In particular, $\PGL(R)$ is topologically simple. \end{cor}

\begin{proof} Of course, Proposition~\ref{proposition:center} and Theorem~\ref{theorem:tbng} together entail that $\PGL(R)$ has topological bounded normal generation. In particular, $\PGL(R)$ is topologically simple by Remark~\ref{remark:topological.simplicity}. \end{proof}

As indicated already in the introduction, Corollary~\ref{corollary:pgl} has been strengthened in the author's more recent work~\cite[Corollary~1.5]{SchneiderSimple}.

The remainder of this section is dedicated to an improvement of Theorem~\ref{theorem:tbng} for certain metric ultraproducts (Corollary~\ref{corollary:ultraproduct}), which vastly extends~\cite[Theorem~9.4]{GMZ25}. The following elementary observation will be useful.

\begin{lem}\label{lemma:sl.dense} Let $K$ be a field and let $n \in \N_{>0}$. Then \begin{displaymath}
	\sup\nolimits_{a \in \GL_{n}(K)} \inf\nolimits_{b \in \SL_{n}(K)} d_{\M_{n}(K)}(a,b) \, \leq \, \tfrac{1}{n} .
\end{displaymath} \end{lem}

\begin{proof} Multiplying one row (or one column) of $a \in \GL_{n}(K)$ by $\det (a)^{-1}$ produces an element $b \in \SL_{n}(K)$ such that $d_{\M_{n}(K)}(a,b) = \rk_{\M_{n}(K)}(a-b) \leq \tfrac{1}{n}$. \end{proof}

\begin{lem}\label{lemma:sl.dense.2} Let $\mathcal{F}$ be an ultrafilter on a set $I$. Consider a family of fields $(K_{i})_{i \in I}$ and a family of positive integers $(n_{i})_{i \in I}$ with $\lim_{i \to \mathcal{F}} n_{i} = \infty$. Then the natural isometric group embedding \begin{align*}
	\prod\nolimits_{i \to \mathcal{F}} (\SL_{n_{i}}(K_{i}),d_{\M_{n_{i}}(K_{i})}) \, &\longrightarrow \, \prod\nolimits_{i \to \mathcal{F}} (\GL_{n_{i}}(K_{i}),d_{\M_{n_{i}}(K_{i})}), \\
	[(a_{i})_{i \in I}] \, &\longmapsto \, [(a_{i})_{i \in I}]
\end{align*} is surjective, i.e., an isometric group isomorphism. \end{lem}

\begin{proof} This is a direct consequence of Lemma~\ref{lemma:sl.dense} and Lemma~\ref{lemma:metric.ultraproduct}\ref{lemma:metric.ultraproduct.1}. \end{proof}

\begin{cor}\label{corollary:ultraproduct} Let $\mathcal{F}$ be a countably incomplete ultrafilter on a set $I$. Suppose that \begin{itemize}
	\item[---\,] either $(R_{i})_{i \in I}$ is a family of non-discrete irreducible, continuous rings,
	\item[---\,] or $(R_{i})_{i \in I} = (\M_{n_{i}}(K_{i}))_{i \in I}$ for some family of fields $(K_{i})_{i \in I}$ and some family of positive integers $(n_{i})_{i \in I}$ with $\lim_{i \to \mathcal{F}} n_{i} = \infty$.
\end{itemize} Let $R$ be the metric ultraproduct of $(R_{i},d_{R_{i}})_{i \in I}$ along $\mathcal{F}$. Then the following hold. \begin{enumerate}
	\item\label{corollary:ultraproduct.1} If $m \in \N$ and $g_{1},\ldots,g_{m} \in \GL(R)$ are such that $\sum\nolimits_{j=1}^{m} d_{R}(g_{j},\cent(R)) > 12$, then \begin{displaymath}
			\qquad \GL(R) \, = \, \Cl_{\GL(R)}(g_{1})\cdots \Cl_{\GL(R)}(g_{m}) .
		\end{displaymath}
	\item\label{corollary:ultraproduct.2} If $g \in \GL(R)\setminus \cent(R)$ and $m \in \N$ are such that $m > \tfrac{12}{d_{R}(g,\cent(R))}$, then
		\begin{displaymath}
			\qquad \GL(R) \, = \, \Cl_{\GL(R)}(g)^{m} .
		\end{displaymath}
	\item\label{corollary:ultraproduct.3} $\PGL(R)$ has bounded normal generation, thus is simple.
\end{enumerate} \end{cor}

\begin{proof} \ref{corollary:ultraproduct.1} Let $m \in \N$ and $g_{1},\ldots,g_{m} \in \GL(R)$ with $\sum\nolimits\nolimits_{j=1}^{m} d_{R}(g_{j},\cent(R)) > 12$. By Proposition~\ref{proposition:ultraproduct} (resp., Lemma~\ref{lemma:sl.dense.2}), there exist $(g_{1,i})_{i \in I},\ldots,(g_{m,i})_{i \in I} \in \prod_{i \in I} \GL(R_{i})$ (resp., $(g_{1,i})_{i \in I},\ldots,(g_{m,i})_{i \in I} \in \prod_{i \in I} \SL_{n_{i}}(K_{i})$) such that \begin{displaymath}
	g_{1} = [(g_{1,i})_{i \in I}], \ \ldots, \ g_{m} = [(g_{m,i})_{i \in I}] .
\end{displaymath} Note that \begin{align*}
	\lim\nolimits_{i\to \mathcal{F}} \sum\nolimits_{j=1}^{m} d_{R_{i}}(g_{j,i},\cent(R_{i})) \, &= \, \sum\nolimits_{j=1}^{m} \lim\nolimits_{i\to \mathcal{F}} d_{R_{i}}(g_{j,i},\cent(R_{i})) \\
		& \stackrel{\ref{proposition:ultraproduct}}{\geq} \, \sum\nolimits_{j=1}^{m} d_{R}(g_{j},\cent(R)) \, > \, 12
\end{align*} and therefore \begin{displaymath}
	I_{0} \, \defeq \, \!\left\{ i \in I \left\vert \, \sum\nolimits_{j=1}^{m} d_{R_{i}}(g_{j,i},\cent(R_{i})) > 12 \right\}\!\right. \, \in \, \mathcal{F} .
\end{displaymath} Since $\mathcal{F}$ is countably incomplete, there exists a sequence $I_{n} \in \mathcal{F}$ $(n \in \N_{>0})$ such that $I_{n+1} \subseteq I_{n}$ for each $n \in \N$ and  $\bigcap_{n \in \N} I_{n} = \emptyset$. Let \begin{displaymath}
	n(i) \, \defeq \, \max \{ n \in \N \mid i \in I_{n} \} \qquad (i \in I_{0}) .
\end{displaymath} Now, consider any $h \in \GL(R)$. We will prove that $h \in \Cl_{\GL(R)}(g_{1})\cdots \Cl_{\GL(R)}(g_{m})$. According to Proposition~\ref{proposition:ultraproduct} (resp., Lemma~\ref{lemma:sl.dense.2}), there exists $(h_{i})_{i \in I} \in \prod_{i \in I} \GL(R_{i})$ (resp., $(h_{i})_{i \in I} \in \prod_{i \in I} \SL_{n_{i}}(K_{i})$) such that $h = [(h_{i})_{i \in I}]$. Applying Theorem~\ref{theorem:tbng} (resp., Corollary~\ref{corollary:index}), for each $i \in I_{0}$, we find $\tilde{g}_{1,i} \in \Cl_{\GL(R_{i})}(g_{1,i}),\ldots,\tilde{g}_{m,i} \in \Cl_{\GL(R_{i})}(g_{m,i})$ such that \begin{displaymath}
	d_{R_{i}}(h_{i},\tilde{g}_{1,i}\cdots \tilde{g}_{m,i}) \, < \, \tfrac{1}{n(i)+1} .
\end{displaymath} Moreover, let $\tilde{g}_{1,i} \defeq 1_{R_{i}}, \ldots , \tilde{g}_{m,i} \defeq 1_{R_{i}}$ for all $i \in I\setminus I_{0}$. Since $I_{0} \in \mathcal{F}$, \begin{displaymath}
	\tilde{g}_{j} \, \defeq \, [(\tilde{g}_{j,i})_{i \in I}] \, \in \, \Cl_{\GL(R)}(g_{j})
\end{displaymath} for each $j \in \{ 1,\ldots,m\}$. Finally, as $I_{n} \in \mathcal{F}$ for every $n \in \N$, \begin{displaymath}
	d_{R}(h,\tilde{g}_{1}\cdots \tilde{g}_{m}) \, \stackrel{\ref{proposition:ultraproduct}}{=} \, \lim\nolimits_{i \to \mathcal{F}} d_{R_{i}}(h_{i},\tilde{g}_{1,i}\cdots \tilde{g}_{m,i}) \, = \, 0 
\end{displaymath} and thus $h = \tilde{g}_{1}\cdots \tilde{g}_{m} \in \Cl_{\GL(R)}(g_{1})\cdots \Cl_{\GL(R)}(g_{m})$, as desired.
	
\ref{corollary:ultraproduct.2} This is an immediate consequence of~\ref{corollary:ultraproduct.1}.
	
\ref{corollary:ultraproduct.3} This follows from~\ref{corollary:ultraproduct.2} and Remark~\ref{remark:topological.simplicity}. \end{proof}

The reader is referred to~\cite[Corollary~1.5]{SchneiderSimple} for a generalization of Corollary~\ref{corollary:ultraproduct}\ref{corollary:ultraproduct.3}.

\section{Geodesic triangularization}\label{section:geodesic.width}

In this section, we establish connectedness properties of unit groups of non-discrete irreducible, continuous rings. Our method pertains to the geometry of algebraic elements of such rings. To be completely explicit, recall that, if $K$ is a field, $R$ is a unital $K$-algebra, and $a \in R$, then the induced evaluation map \begin{displaymath}
	K[X] \, \longrightarrow \, R, \quad p = \sum\nolimits_{i=0}^{n} c_{i}X^{i} \, \longmapsto \, p(a) \defeq p_{R}(a) \defeq \sum\nolimits_{i=0}^{n} c_{i}a^{i}
\end{displaymath} is a unital $K$-algebra homomorphism.

\begin{definition} Let $K$ be a field and let $R$ be a unital $K$-algebra. For $S \subseteq K[X]$, we let $\alg_{R}(S) \defeq \{ a \in R \mid \forall p \in S \colon \, p(a) = 0 \}$. The members of \begin{displaymath}
	\A(R) \, \defeq \, \bigcup\nolimits_{p \in K[X]\setminus \{ 0\}}\alg_{R}(\{ p \})
\end{displaymath} are referred to as \emph{algebraic} or \emph{$K$-algebraic} elements of $R$. \end{definition}

Our approach is based on \emph{continuous triangularization} as developed in~\cite{SchneiderGAFA}, thus involving the following notions and objects.

\begin{definition}[{\cite[Definition~7.14]{SchneiderGAFA}}] Let $R$ be a unital ring. A \emph{nest} in $R$ is a chain in $(\E(R),{\leq})$. Let $\Nest(R)$ denote the set of all nests in $R$, and let $\Nestmax(R)$ denote the set of all maximal elements in the partially ordered set $(\Nest(R),{\subseteq})$. \end{definition}

\begin{lem}[\cite{SchneiderGAFA}]\label{lemma:homomorphism} Let $R$ be a unital ring and let $E \subseteq \E(R)$. Then \begin{align*}
	R_{E} \, &\defeq \, \{ a \in R \mid \forall e \in E \colon \, eae = ae \} \\
	&= \, \{ a \in R \mid \forall e \in E \colon \, (1-e)a(1-e) = (1-e)a \}
\end{align*} is a unital $\cent(R)$-subalgebra of $R$. Moreover, if $e \in \{ f-f' \mid f,f' \in E \cup \{ 0,1\}, \, f' \leq f \}$, then  \begin{displaymath}
	R_{E} \, \longrightarrow \, eRe, \quad a \, \longmapsto \, eae
\end{displaymath} is a unital $\cent(R)$-algebra homomorphism. \end{lem}

\begin{proof} The first part follows from~\cite[Lemma~9.1]{SchneiderGAFA} and~\cite[Proposition~9.3(2)]{SchneiderGAFA}, while the second is proved in~\cite[Lemma~10.5(1)]{SchneiderGAFA}. \end{proof}

The \emph{continuous triangularization theorem}~\cite[Corollary 9.12(1)]{SchneiderGAFA} (resp., its invertible version~\cite[Corollary 9.12(2)]{SchneiderGAFA}) asserts that every algebraic element of an irreducible, continuous ring $R$ (resp., of its unit group $\GL(R)$) is contained in $R_{E}$ (resp., $\GL(R_{E})$) for some $E \in \Nestmax(R)$. Let us briefly summarize some direct topological ramifications of the results of~\cite{SchneiderGAFA}.

\begin{remark}\label{remark:summary} Let $R$ be a non-discrete irreducible, continuous ring. \begin{enumerate}
	\item\label{remark:summary.1} If $E \in \Nestmax (R)$, then $\GL(R_{E})$ is contractible, in particular connected, with respect to the rank topology by~\cite[Lemma~11.1]{SchneiderGAFA} and~\cite[Lemma~11.4]{SchneiderGAFA}.
	\item\label{remark:summary.2} Since \begin{displaymath}
			\qquad \GL(R) \, = \, \overline{\bigcup\nolimits_{E \in \Nestmax(R)} \GL(R_{E})}
	\end{displaymath} by~\cite[Corollary~9.13]{SchneiderGAFA} and $1 \in \bigcap\nolimits_{E \in \Nestmax(R)} \GL(R_{E})$, it follows from~\ref{remark:summary.1} that $\GL(R)$ is connected with respect to the rank topology (see, e.g.,~\cite[I, \S11.1, Proposition~2, p.~108]{bourbaki} and~\cite[I, \S11.1, Proposition~1, p.~108]{bourbaki}).
\end{enumerate}\end{remark}

We proceed to this section's technical core.

\begin{lem}\label{lemma:geodesic} Let $R$ be a non-discrete irreducible, continuous ring, let $a,b \in R$ and $e \in \E(R)$ with $(b-a)R = eR$, let $E \in \Nestmax(eRe)$ and \begin{displaymath}
	e_{t} \, \defeq \, ({\rk_{R}}\vert_{E})^{-1}(t) \qquad (t \in [0,\rk_{R}(e)]) .
\end{displaymath} Then \begin{displaymath}
	\gamma \colon \, [0,\rk_{R}(e)] \, \longrightarrow \, R, \quad t \, \longmapsto \, e_{t}b + (1-e_{t})a
\end{displaymath} is a geodesic in $(R,d_{R})$ from $a$ to $b$. Suppose that $a \in R_{E}$ and $b \in R_{\{ 1-f \mid f \in E \}}$. Then, moreover, the following hold. \begin{enumerate}
	\item\label{lemma:geodesic.1} If $a,b \in \GL(R)$, then $\im(\gamma) \subseteq \GL(R)$.
	\item\label{lemma:geodesic.2} If $S \subseteq \cent(R)[X]$ and $a,b \in \alg_{R}(S)$, then $\im(\gamma) \subseteq \alg_{R}(S)$.
\end{enumerate} \end{lem}

\begin{proof} From Lemma~\ref{lemma:corner.ring}\ref{lemma:corner.ring.3} and~\cite[Corollary~7.20(2)]{SchneiderGAFA}, we know that \begin{displaymath}
	([0,\rk_{R}(e)],{\leq}) \, \longrightarrow \, (E,{\leq}), \quad t \, \longmapsto \, e_{t}
\end{displaymath} is a well-defined order isomorphism. Evidently, $\gamma(0) = a$ and \begin{align*}
	\gamma(\rk_{R}(e)) \, &= \, eb+(1-e)a \, = \, eb+(1-e)a + (1-e)e(b-a) \\
		& \stackrel{b-a \in eR}{=} \, eb+(1-e)a + (1-e)(b-a) \, = \, eb + (1-e)b \, = \, b .
\end{align*} Moreover, if $s,t \in [0,\rk_{R}(e)]$ and $s \leq t$, then \begin{align*}
	d_{R}(\gamma(s),\gamma(t)) \, &= \, \rk_{R}(\gamma(s)-\gamma(t)) \, = \, \rk_{R}((e_{t}-e_{s})b-(e_{t}-e_{s})a) \\
		& = \, \rk_{R}((e_{t}-e_{s})(b-a)) \, \stackrel{(b-a)R = eR}{=} \, \rk_{R}((e_{t}-e_{s})e) \\
		& \stackrel{e_{s},e_{t}\leq e}{=} \, \rk_{R}(e_{t}-e_{s}) \, \stackrel{\ref{lemma:pseudo.rank.function}\ref{lemma:pseudo.rank.function.1}}{=} \, \rk_{R}(e_{t}) - \rk_{R}(e_{s}) \, = \, t-s \, = \, \vert s-t \vert .
\end{align*} This shows that $\gamma$ is an isometric embedding, thus constitutes a geodesic from $a$ to $b$ in $(R,d_{R})$, as desired. 

Turning to the proofs of~\ref{lemma:geodesic.1} and~\ref{lemma:geodesic.2}, let us consider $E' \defeq \{ 1-f \mid f \in E \} \subseteq \E(R)$. For each $t \in [0,\rk_{R}(e)]$, since $e_{t}Re_{t} \times (1-e_{t})R(1-e_{t}) \to R, \, (y,x) \mapsto y+x$ is a unital $\cent(R)$-algebra homomorphism according to~\cite[Remark~10.1(1)]{SchneiderGAFA}, we conclude using Lemma~\ref{lemma:homomorphism} that \begin{align*}
	\phi_{t} \colon \, R_{E} \times R_{E'} \, \longrightarrow \, R, \quad (x,y) \, \longmapsto \, e_{t}ye_{t} + (1-e_{t})x(1-e_{t}) \stackrel{\ref{lemma:homomorphism}}{=} e_{t}y + (1-e_{t})x
\end{align*} is a unital $\cent(R)$-algebra homomorphism, too. Now, suppose that $a \in R_{E}$ and $b \in R_{E'}$. Then $\phi_{t}(a,b) = \gamma(t)$. We continue with separate arguments.

\ref{lemma:geodesic.1} If $a,b \in \GL(R)$, thus $(a,b) \in \GL(R_{E}) \times \GL(R_{E'})$ by~\cite[Proposition~9.5]{SchneiderGAFA}, then \begin{displaymath}
	\gamma(t) \, = \, \phi_{t}(a,b) \, \in \, \phi_{t}(\GL(R_{E}) \times \GL(R_{E'})) \, = \, \phi_{t}(\GL(R_{E} \times R_{E'})) \, \subseteq \, \GL(R) 
\end{displaymath} for every $t \in [0,\rk_{R}(e)]$.

\ref{lemma:geodesic.2} Let $S \subseteq K[X]$ with $a,b \in \alg_{R}(S)$. Now, if $t \in [0,\rk_{R}(e)]$, then \begin{displaymath}
	p(\gamma(t)) \, = \, p_{R}(\phi_{t}(a,b)) \, = \, \phi_{t}(p_{R_{E}\times R_{E'}}(a,b)) \, = \, \phi_{t}(p_{R_{E}}(a),p_{R_{E'}}(b)) \, = \, \phi_{t}(0,0) \, = \, 0
\end{displaymath} for each $p \in S$ and thus $\gamma(t) \in \alg_{R}(S)$. \end{proof}

\begin{prop}\label{proposition:star} Let $R$ be a non-discrete irreducible, continuous ring, $S \subseteq \cent (R)[X]$ with $\emptyset \ne S \ne \{ 0 \}$, and $c \in \alg_{\cent (R)}(S)$. Then $(\alg_{R}(S),d_{R},c)$ is star-shaped. \end{prop}

\begin{proof} Let $a \in \alg_{R}(S)\setminus \{ c \}$. Due to Remark~\ref{remark:regular}, we find $e \in \E(R)$ with $(c-a)R = eR$. Note that $e \ne 0$ as $c \ne a$. Therefore, Lemma~\ref{lemma:corner.ring}\ref{lemma:corner.ring.3} asserts that $eRe$ is a non-discrete irreducible, continuous ring with $\cent (eRe) = \cent(R)e \cong \cent(R)$. Moreover, \begin{displaymath}
	aeR \, = \, a(c-a)R \, = \, (c-a)aR \, \subseteq \, (c-a)R \, = \, eR
\end{displaymath} and hence $eae = ae$, i.e., $a \in R_{\{ e \}}$. Since $R_{\{ e \}} \to eRe, \, x \mapsto exe = xe$ is a unital $\cent(R)$-algebra homomorphism by Lemma~\ref{lemma:homomorphism}, we see that $p_{eRe}(ae) = p_{R_{\{ e \}}}(a)e = p_{R}(a)e = 0$ for each $p \in S$. By assumption, $S$ contains at least one non-zero element, so $ae$ is algebraic. Due to~\cite[Corollary 9.12(1)]{SchneiderGAFA}, thus there exists $E \in \Nestmax(eRe)$ such that $ae \in (eRe)_{E} = eRe \cap R_{E}$. For each $f \in E$, we infer that \begin{displaymath}
	af \, \stackrel{f \leq e}{=} \, aef \, \stackrel{ae \in R_{E}}{=} \, faef \, \stackrel{f \leq e}{=} \, faf .
\end{displaymath} That is, $a \in R_{E}$. Also, $c \in \cent (R) \subseteq R_{\{ 1-f \mid f \in E\}}$. Hence, by Lemma~\ref{lemma:geodesic}, there exists a geodesic from $a$ to $c$ in $(\alg_{R}(S),d_{R})$. \end{proof}

\begin{cor} Let $R$ be a non-discrete irreducible, continuous ring. If $K \defeq \cent (R)$ is algebraically closed, then the metric space $(\alg_{R}(I),d_{R})$ is star-shaped for every proper, non-zero ideal $I \leq K[X]$. \end{cor}

\begin{proof} Consider any proper, non-zero ideal $I \leq K[X]$. Since $K[X]$ is a principal ideal domain, we find $p \in K[X]\setminus K$ such that $I = K[X]p$. Now, if $K$ is algebraically closed, then $\alg_{R}(I) = \alg_{R}(\{ p \})$ has a non-empty intersection with $K$, whence $(\alg_{R}(I),d_{R})$ is star-shaped due to Proposition~\ref{proposition:star}. \end{proof}

Given a unital ring $R$, we consider $\I(R) \defeq \left\{ a \in R \left\vert \, a^{2} = 1 \right\} \! \right. \subseteq \GL(R)$.

\begin{cor}\label{corollary:star} Let $R$ be a non-discrete irreducible, continuous ring. Then the pointed metric space $(\I(R),d_{R},1)$ is star-shaped. In particular, $\I(R) \subseteq \Delta(\GL(R),{d_{R}})$. \end{cor}

\begin{proof} This follows upon applying Proposition~\ref{proposition:star} to $\alg_{R}(\{ X^{2}-1 \}) = \I(R)$. \end{proof}

Similarly, Lemma~\ref{lemma:geodesic} constitutes a central ingredient in the following insight.

\begin{thm}\label{theorem:geodesic} Let $R$ be a non-discrete irreducible, continuous ring. Then \begin{displaymath}
	\GL(R) \cap \A(R) \, \subseteq \, \Delta(\GL(R),{d_{R}}) .
\end{displaymath} In particular, \begin{displaymath}
	\GL(R) \, = \, \overline{\Delta(\GL(R),{d_{R}})} .
\end{displaymath} \end{thm}

\begin{proof} Let $a \in \GL(R) \cap \A(R)$. If $a=1$, then clearly $a \in \Delta(\GL(R),{d_{R}})$. Henceforth, we will assume that $a \ne 1$. By Remark~\ref{remark:regular}, there exists $e \in \E(R)$ with $(1-a)R = eR$. Notice that $e \ne 0$ as $a \ne 1$. So, Lemma~\ref{lemma:corner.ring}\ref{lemma:corner.ring.3} asserts that $eRe$ is a non-discrete irreducible, continuous ring with $\cent (eRe) = \cent(R)e \cong \cent(R)$. Furthermore, \begin{displaymath}
	aeR \, = \, a(1-a)R \, = \, (1-a)aR \, \subseteq \, (1-a)R \, = \, eR
\end{displaymath} and hence $eae = ae$, i.e., $a \in R_{\{ e\}}$. As $a \in \A(R)$, there exists $p \in \cent(R)[X]\setminus \{ 0 \}$ such that $p(a)=0$. Since $R_{\{ e \}} \to eRe, \, x \mapsto exe = xe$ is a unital $\cent(R)$-algebra homomorphism by Lemma~\ref{lemma:homomorphism}, we deduce that $p_{eRe}(ae) = p_{R_{\{ e \}}}(a)e = p_{R}(a)e = 0$, so $ae$ is algebraic. According to~\cite[Corollary 9.12(1)]{SchneiderGAFA}, we find $E \in \Nestmax(eRe)$ such that $ae \in (eRe)_{E} = eRe \cap R_{E}$. For every $f \in E$, it follows that \begin{displaymath}
	af \, \stackrel{f \leq e}{=} \, aef \, \stackrel{ae \in R_{E}}{=} \, faef \, \stackrel{f \leq e}{=} \, faf .
\end{displaymath} Thus, $a \in R_{E}$. Of course, $1 \in R_{\{ 1-f \mid f \in E\}}$. Therefore, thanks to Lemma~\ref{lemma:geodesic}, we find a geodesic in $(\GL(R),d_{R})$ from $a$ to $1$, which shows that $a \in \Delta(\GL(R),{d_{R}})$. This proves the first assertion. Since $\GL(R) \cap \A(R)$ is dense in $(\GL(R),d_{R})$ according to~\cite[Proposition~8.5]{SchneiderGAFA}, the second assertion follows at once. \end{proof}

\begin{cor}\label{corollary:length.space} Let $R$ be a non-discrete irreducible, continuous ring. Then the metric space $(\GL(R),d_{R})$ is a length space. In particular, $\GL(R)$ is both path-connected and locally path-connected with respect to the rank topology. \end{cor}

\begin{proof} From Remark~\ref{remark:properties.pseudo.rank.function}\ref{remark:unit.group}, we know that $(\GL(R),d_{R})$ is a metric group. Moreover, $(\GL(R),d_{R})$ is complete by Theorem~\ref{theorem:unique.rank.function}\ref{theorem:unique.rank.function.1} and Remark~\ref{remark:properties.pseudo.rank.function}\ref{remark:unit.group}. Hence, $(\GL(R),d_{R})$ is a length space due to Theorem~\ref{theorem:geodesic} and Lemma~\ref{lemma:approximate.midpoints}. It follows that $\GL(R)$ is both path-connected and locally path-connected with respect to the rank topology: while the former is trivial, the latter is a consequence of Lemma~\ref{lemma:length.space.locally.path.connected}. \end{proof}

\begin{cor} Let $\mathcal{F}$ be a countably incomplete ultrafilter on a set $I$ and let $(R_{i})_{i \in I}$ be a family of irreducible, continuous rings such that $R \defeq \prod_{i \to \mathcal{F}} R_{i}$ is non-discrete. Then $(\GL(R),{d_{R}})$ is geodesic. \end{cor}

\begin{proof} Thanks to Corollary~\ref{corollary:length.space}, we know that $(\GL(R),{d_{R}})$ is a length space. Since $(\GL(R),{d_{R}}) \cong \prod_{i \to \mathcal{F}} (\GL(R_{i}),{d_{R_{i}}})$ by Proposition~\ref{proposition:ultraproduct}, it follows that $(\GL(R),{d_{R}})$ is geodesic due to Proposition~\ref{proposition:length.spaces}\ref{proposition:length.spaces.1} and Lemma~\ref{lemma:metric.ultraproduct}\ref{lemma:metric.ultraproduct.3}.  \end{proof}

We conclude this section by deducing a uniform finite upper bound on the geodesic width of unit groups of non-discrete irreducible, continuous rings (Theorem~\ref{theorem:geodesic.width}). This result will become a direct consequence of some previous work (Theorem~\ref{theorem:locally.algebraic}), once the conclusion of Theorem~\ref{theorem:geodesic} has been extended to a larger class of elements (Lemma~\ref{lemma:locally.algebraic}) using the following standard decomposition technique. Let us recall that the set $\Pfin(I)$ of all finite subsets of a set $I$, equipped with the partial order defined by inclusion, constitutes a directed set.

\begin{lem}\label{lemma:convergence.sequences} Let $\rho$ be a rank function on a regular ring $R$ such that $(R,{d_{\rho}})$ is complete. Let $E$ be a set of pairwise orthogonal elements of $\E(R)$. Then \begin{displaymath}
	\prod\nolimits_{e \in E} eRe \, \longrightarrow \, R, \quad a \, \longmapsto \, \sum\nolimits_{e \in E} a_{e} \defeq \lim\nolimits_{F \to \Pfin(E)} \sum\nolimits_{e \in F} a_{e}
\end{displaymath} is a well-defined ring embedding and \begin{displaymath}
	\forall a \in \prod\nolimits_{e \in E} eRe \colon \quad \rho\!\left( \sum\nolimits_{e \in E} a_{e} \right)\! \, = \, \sum\nolimits_{e \in E} \rho(a_{e}) .
\end{displaymath} In particular, if $\sum\nolimits_{e \in E} e = 1$, then \begin{displaymath}
	\prod\nolimits_{e \in E} \GL(eRe) \, = \, \! \left\{ a \in \prod\nolimits_{e \in E} eRe \left\vert \, \sum\nolimits_{e \in E} a_{e} \in \GL(R) \right\} . \right.
\end{displaymath} \end{lem}

\begin{proof} The first assertion is a consequence of~\cite[Lemma~5.6]{BernardSchneider}. Combining the continuity of the map $(R,d_{\rho}) \to [0,1], \, a \mapsto \rho(a) = d_{\rho}(a,0)$ with~\cite[Lemma~5.1(B)]{BernardSchneider}, we see that \begin{align*}
	\rho\!\left( \sum\nolimits_{e \in E} a_{e} \right)\! \, &= \, \lim\nolimits_{F \to \Pfin(E)} \rho\!\left( \sum\nolimits_{e \in F} a_{e} \right) \\
	& = \,  \lim\nolimits_{F \to \Pfin(E)} \sum\nolimits_{e \in F} \rho(a_{e}) \, = \, \sum\nolimits_{e \in E} \rho(a_{e}) 
\end{align*} for all $a \in \prod\nolimits_{e \in E} eRe$. Finally, if $\sum\nolimits_{e \in E} e = 1$ and thus $\sum\nolimits_{e \in E} \rho(e) = \rho\!\left( \sum\nolimits_{e \in E} e \right)\! = 1$, then we conclude that \begin{align*}
	\prod\nolimits_{e \in E} \GL(eRe) \, &\stackrel{\ref{lemma:corner.ring}\ref{lemma:corner.ring.2}+\ref{remark:properties.pseudo.rank.function}\ref{remark:unit.group}}{=} \, \! \left. \left\{ a \in \prod\nolimits_{e \in E} eRe \, \right\vert \forall e \in E \colon \, \rho(a_{e}) = \rho(e) \right\} \\
	&= \, \! \left. \left\{ a \in \prod\nolimits_{e \in E} eRe \, \right\vert \sum\nolimits_{e \in E} \rho(a_{e}) = 1 \right\} \\
	&= \, \! \left. \left\{ a \in \prod\nolimits_{e \in E} eRe \, \right\vert \rho\!\left( \sum\nolimits_{e \in E} a_{e} \right)\! =  1 \right\} \\
	&\stackrel{\ref{remark:properties.pseudo.rank.function}\ref{remark:unit.group}}{=} \, \! \left\{ a \in \prod\nolimits_{e \in E} eRe \left\vert \, \sum\nolimits_{e \in E} a_{e} \in \GL(R) \right\} \right. . \qedhere
\end{align*} \end{proof}

\begin{definition}\label{definition:locally.algebraic} Let $R$ be an irreducible, continuous ring. An element $a \in R$ is called \emph{locally algebraic} if there exist a (necessarily countable\footnote{cf.~\cite[Remark~5.5]{BernardSchneider}}) set $E$ of pairwise orthogonal elements of $\E(R)\setminus \{ 0 \}$ and $(a_{e})_{e \in E} \in \prod_{e \in E} \A(eRe)$ such that $a = \sum\nolimits_{e \in E} a_{e}$. \end{definition}

\begin{lem}\label{lemma:locally.algebraic} Let $R$ be a non-discrete irreducible, continuous ring. If $a \in \GL(R)$ is locally algebraic, then $a \in \Delta(\GL(R),{d_{R}})$. \end{lem}

\begin{proof} Let $a \in \GL(R)$ be locally algebraic. Then we find a (necessarily countable) set $E$ of pairwise orthogonal elements of $\E(R)\setminus \{ 0 \}$ and $b \in \prod_{e \in E} \A(eRe)$ such that the map \begin{displaymath}
	\phi \colon \, \prod\nolimits_{e \in E} eRe \, \longrightarrow \, R, \quad x \, \longmapsto \, \sum\nolimits_{e \in E} x_{e}
\end{displaymath} satisfies $\phi(b) = a$. Since $\phi$ is a ring homomorphism by Lemma~\ref{lemma:convergence.sequences}, we see that $\sum_{e \in E} e = \phi((e)_{e \in E}) \in \E(R)$, while also \begin{displaymath}
	\left( \sum\nolimits_{e \in E} e\right) \! a \! \left( \sum\nolimits_{e \in E} e\right) \, = \, \phi\!\left( (e)_{e \in E} \, b \, (e)_{e \in E}\right) \, = \, \phi(b) \, = \, a \, \in \, \GL(R) 
\end{displaymath} and thus $\sum\nolimits_{e \in E} e \in \GL(R)$. As $\E(R) \cap \GL(R) = \{ 1 \}$, it follows that $\sum_{e \in E} e = 1$. For each $e \in E$, we conclude that $b_{e} \in \GL(eRe)$ by Lemma~\ref{lemma:convergence.sequences}, and then apply Lemma~\ref{lemma:corner.ring}\ref{lemma:corner.ring.3} and Theorem~\ref{theorem:geodesic} to find a geodesic from $b_{e}$ to $e$ in $(\GL(eRe),d_{eRe})$, which in turn constitutes a geodesic from $b_{e}$ to $e$ in $(\GL(eRe),d_{R})$. Therefore, according to Remark~\ref{remark:geodesics.in.products}, there exists a geodesic $\gamma$ from $b$ to $(e)_{e \in E}$ in $G \defeq \prod\nolimits_{e \in E}\GL(eRe)$ with respect to the metric \begin{displaymath}
	G \times G \, \longrightarrow \, [0,1], \quad (x,y) \, \longmapsto \, \sum\nolimits_{e \in E} d_{R}(x_{e},y_{e}) .
\end{displaymath} Since $\phi(G) \subseteq \GL(R)$ by Lemma~\ref{lemma:convergence.sequences} and \begin{align*}
	d_{R}(\phi(x),\phi(y)) \, &= \, \rk_{R}(\phi(x)-\phi(y)) \, \stackrel{\ref{lemma:convergence.sequences}}{=} \, \rk_{R}(\phi(x-y)) \, = \, \rk_{R}\!\left( \sum\nolimits_{e \in E} x_{e}-y_{e} \right) \\
	&\stackrel{\ref{theorem:unique.rank.function}+\ref{lemma:convergence.sequences}}{=} \, \sum\nolimits_{e \in E} \rk_{R}(x_{e}-y_{e}) \, = \, \sum\nolimits_{e \in E} d_{R}(x_{e},y_{e}) 
\end{align*} for all $x,y \in G$, we see that $\phi \circ \gamma$ constitutes a geodesic in $(\GL(R),d_{R})$ from $\phi(b) = a$ to $\phi((e)_{e \in E}) = 1$. Consequently, $a \in \Delta(\GL(R),{d_{R}})$. \end{proof}

\begin{thm}[\cite{BernardSchneider}]\label{theorem:locally.algebraic} If $R$ is a non-discrete irreducible, continuous ring, then every element of $\GL(R)$ is a product of $7$ locally algebraic elements of $\GL(R)$. \end{thm}

\begin{proof} This follows from~\cite[Theorem~1.1]{BernardSchneider} and~\cite[Theorem~8.9]{BernardSchneider}. \end{proof}

We arrive at this section's final result.

\begin{thm}\label{theorem:geodesic.width} If $R$ is a non-discrete irreducible, continuous ring, then \begin{displaymath}
	\GL(R) \, = \, \Delta (\GL(R),{d_{R}})^{7} .
\end{displaymath} \end{thm}

\begin{proof} This follows from Theorem~\ref{theorem:locally.algebraic} and Lemma~\ref{lemma:locally.algebraic}. \end{proof}

\section{Escape dynamics}\label{section:escape.dynamics}

This section is dedicated to the proof of Theorem~\ref{theorem:bounded} and its ramifications regarding escape dynamics (Corollary~\ref{corollary:no.escape}), homomorphism rigidity (Corollaries~\ref{corollary:homomorphism.rigidity} and~\ref{corollary:rigidity}), Bourbaki boundedness (Corollary~\ref{corollary:bounded}), and representation theory (Corollary~\ref{corollary:strongly.exotic}). The following construction of small subgroups will be useful.

\begin{lem}[cf.~{\cite[Lemma~5.2]{BernardSchneider}}]\label{lemma:gamma} Let $R$ be a unital ring and let $e \in \E(R)$. Then \begin{displaymath}
	\Gamma_{R}(e) \, \defeq \, \GL(eRe) + 1-e \, = \, \GL(R) \cap (eRe+1-e)
\end{displaymath} is a subgroup of $\GL(R)$ and \begin{displaymath}
	\phi \colon \, \GL(eRe) \, \longrightarrow \, \Gamma_{R}(e), \quad a \, \longmapsto \, a+1-e
\end{displaymath} is an isomorphism. If $R$ is regular and $\rho$ is a pseudo-rank function on $R$, then $d_{\rho}(\phi(a),\phi(b)) = d_{\rho}(a,b)$ for all $a,b \in \GL(eRe)$. \end{lem}

\begin{proof} The first part is due to~\cite[Lemma~5.2]{BernardSchneider}. Moreover, if $R$ is regular and $\rho$ is a pseudo-rank function on $R$, then, for all $a,b \in \GL(eRe)$, \begin{displaymath}
	d_{\rho}(\phi(a),\phi(b)) \, = \, \rho(\phi(a)-\phi(b)) \, = \, \rho(a-b) \, = \, d_{\rho}(a,b) . \qedhere
\end{displaymath} \end{proof}

\begin{lem}\label{lemma:small.subgroup} Let $R$ be a unital ring, let $e,f \in \E(R)$ with $e \perp f$. \begin{enumerate}
	\item\label{lemma:small.subgroup.action} The map \begin{displaymath}
				\qquad \alpha \colon \, eRf \times \Gamma_{R}(f)\, \longrightarrow \, eRf, \quad (x,a) \, \longmapsto \, xa
			\end{displaymath} is a right action of $\Gamma_{R}(f)$ on $eRf$ by left $eRe$-module automorphisms.
	\item\label{lemma:small.subgroup.isomorphism} $\Gamma_{R}(f) + eRf$ is a subgroup of $\GL(R)$ and \begin{displaymath}
				\qquad \Gamma_{R}(f) \ltimes_{\alpha} eRf \, \longrightarrow \, \Gamma_{R}(f) + eRf, \quad (a,x) \, \longmapsto \, a+x
	\end{displaymath} is an isomorphism. 
	\item\label{lemma:small.subgroup.ball} If $R$ is regular and $\rho$ is a pseudo-rank function on $R$, then $\Gamma_{R}(f) + eRf$ is contained in $\ball_{\rho(f)}(\GL(R),{d_{\rho}})$, in particular \begin{displaymath}
				\qquad \Gamma_{R}(f) + eRf \, \subseteq \, \trap\!\left(\ball_{\rho(f)}(\GL(R),{d_{\rho}})\right) .
\end{displaymath} \end{enumerate} \end{lem}

\begin{proof} \ref{lemma:small.subgroup.action} Note that $\alpha$ is well defined: indeed, if $x \in eRf$ and $a \in \Gamma_{R}(f)$, then $xa = exf(faf+1-f) = exfaf \in eRf$. Routine calculations show that $\alpha$ is a group action by left $eRe$-module automorphisms.
	
\ref{lemma:small.subgroup.isomorphism} Consider $G \defeq \Gamma_{R}(f) \ltimes_{\alpha} eRf$ and $\phi \colon G \to R, \, (a,x) \mapsto a+x$. Then $\phi(1,0) = 1$. Moreover, if $(a,x),(b,y) \in G$, then \begin{align*}
	\phi((a,x)\cdot (b,y)) \, &= \, \phi(ab,xb+y) \, = \, ab+xb+y \, = \, ab+xb+ey \\
	&\stackrel{e \perp f}{=} \, ab+xb+(faf+1-f)ey+xfey \\
	&= \, ab+xb+ay+xy \, = \, (a+x)(b+y) \, = \, \phi(a,x)\phi(b,y) .
\end{align*} Since $G$ is a group, it follows that $\Gamma_{R}(f) + eRf = \phi(G)$ is contained in and constitutes a subgroup of $\GL(R)$ and that $G \to \GL(R), \, g \mapsto \phi(g)$ is a homomorphism. Finally, if $(a,x) \in G$ and $\phi(a,x)=1$, then \begin{displaymath}
	x \, = \, ex+e-e \, \stackrel{e \perp f}{=} \, ex+ea-e \, = \, e(a+x)-e \, = \, e\phi(a,x)-e \, = \, e-e \, = \, 0
\end{displaymath} and therefore $a = a+x = \phi(a,x) = 1$, whence $(a,x) = (1,0)$. This shows that $\phi$ is injective, thus induces the desired isomorphism.

\ref{lemma:small.subgroup.ball} Let $R$ be regular and $\rho$ be a pseudo-rank function on $R$. Then, for all $a \in \Gamma_{R}(f)$ and $x \in eRf$, \begin{align*}
	d_{\rho}(a+x,1) \, &= \, \rho(faf + 1-f + exf -1) \\
	& = \, \rho(faf-exf-f) \, = \, \rho((fa-ex-1)f) \, \leq \, \rho(f) .
\end{align*} That is, $\Gamma_{R}(f) + eRf \subseteq \ball_{\rho(f)}(\GL(R),{d_{\rho}})$. As $\Gamma_{R}(f) + eRf$ is a subgroup of $\GL(R)$ due to~\ref{lemma:small.subgroup.isomorphism}, this entails that $\Gamma_{R}(f) + eRf \subseteq \trap\!\left(\ball_{\rho(f)}(\GL(R),{d_{\rho}})\right)$. \end{proof}

\begin{lem}[cf.~{\cite[Lemma~11.1(1)]{SchneiderGAFA}}]\label{lemma:idempotent.difference} Let $R$ be a unital ring. Moreover, let $E \in \Nest (R)$ and $e \in \{ f-f' \mid f,f' \in E \cup \{ 0,1 \}, \, f' \leq f \}$. Then \begin{displaymath}
	\pi_{e} \colon \, \GL(R_{E}) \, \longrightarrow \, \GL(R_{E}) , \quad a \, \longmapsto \, eae + 1-e
\end{displaymath} is a group endomorphism with $\pi_{e}(\GL(R_{E})) \subseteq \Gamma_{R}(e)$. \end{lem}

\begin{proof} By~\cite[Lemma~11.1(1)]{SchneiderGAFA}, $\pi_{e}$ is a well-defined group endomorphism. In turn, \begin{displaymath}
	\pi_{e}(\GL(R_{E})) \, \subseteq \, \GL(R) \cap (eRe + 1-e) \, \stackrel{\ref{lemma:gamma}}{=} \, \Gamma_{R}(e) .\qedhere
\end{displaymath} \end{proof}

\begin{thm}\label{theorem:bounded} Let $R$ be a non-discrete irreducible, continuous ring, and $n \in \N_{>0}$. Then \begin{displaymath}
	\GL(R) \cap \A(R) \, \subseteq \, \trap\!\left(\ball_{1/n}(\GL(R),d_{R})\right)^{n} .
\end{displaymath} In particular, \begin{displaymath}
	\GL(R) \, = \, \overline{\trap\!\left(\ball_{1/n}(\GL(R),d_{R})\right)^{n}} \, = \, \trap\!\left(\ball_{1/n}(\GL(R),d_{R})\right)^{n+1} .
\end{displaymath} \end{thm}

\begin{proof} Let $g \in \GL(R) \cap \A(R)$. By~\cite[Corollary~9.12(2)]{SchneiderGAFA}, we find $E \in \Nestmax (R)$ such that $g \in \GL(R_{E})$. Due to~\cite[Corollary~7.20(2)]{SchneiderGAFA}, the map ${\rk_{R}}\vert_{E} \colon (E,{\leq}) \to ([0,1],{\leq})$ is an order isomorphism. For each $i \in \{ 0,\ldots,n\}$, we let \begin{displaymath}
	e_{i} \, \defeq \, \left( {\rk_{R}}\vert_{E}\right)^{-1}\!\left( \tfrac{i}{n} \right)\! \, \in \, E .
\end{displaymath} Furthermore, for each $i \in \{ 1,\ldots,n\}$, consider $f_{i} \defeq e_{i}-e_{i-1} \in \E(R)$ and note that \begin{displaymath}
	\rk_{R}(f_{i}) \, = \, \rk_{R}(e_{i}-e_{i-1}) \, \stackrel{\ref{lemma:pseudo.rank.function}\ref{lemma:pseudo.rank.function.1}}{=} \, \rk_{R}(e_{i}) - \rk_{R}(e_{i-1}) \, = \, \tfrac{i}{n} - \tfrac{i-1}{n} \, = \, \tfrac{1}{n} .
\end{displaymath} Using Lemma~\ref{lemma:idempotent.difference}, we define $g_{1} \defeq \pi_{e_{1}}(g) \in \GL(R_{E})$ and recursively \begin{displaymath}
	g_{i} \, \defeq \, \pi_{e_{i}}\!\left( gg_{1}^{-1}\cdots g_{i-1}^{-1}\right) \, \in \, \GL(R_{E}) \qquad (i \in \{ 2,\ldots,n \}) .
\end{displaymath} Evidently, \begin{displaymath}
	g_{1} \, = \, \pi_{e_{1}}(g) \, \stackrel{\ref{lemma:idempotent.difference}}{\in} \, \Gamma_{R}(e_{1}) \, \stackrel{\ref{lemma:small.subgroup}\ref{lemma:small.subgroup.ball}}{\subseteq} \, \trap\!\left(\ball_{1/n}(\GL(R),d_{R})\right).
\end{displaymath} Moreover, for each $i \in \{ 2,\ldots,n\}$, \begin{align*}
	gg_{1}^{-1}\cdots g_{i-1}^{-1} \, &= \, gg_{1}^{-1}\cdots g_{i-2}^{-1}\pi_{e_{i-1}}\!\left(gg_{1}^{-1}\cdots g_{i-2}^{-1}\right)^{-1} \\
	& = \, gg_{1}^{-1}\cdots g_{i-2}^{-1}\pi_{e_{i-1}}\!\left(\!\left(gg_{1}^{-1}\cdots g_{i-2}^{-1}\right)^{-1}\right) \\
	&= \, gg_{1}^{-1}\cdots g_{i-2}^{-1}e_{i-1}\!\left(gg_{1}^{-1}\cdots g_{i-2}^{-1}\right)^{-1}\!e_{i-1} + gg_{1}^{-1}\cdots g_{i-2}^{-1}(1-e_{i-1}) \\
	&= \, gg_{1}^{-1}\cdots g_{i-2}^{-1}\!\left(gg_{1}^{-1}\cdots g_{i-2}^{-1}\right)^{-1}\!e_{i-1} + gg_{1}^{-1}\cdots g_{i-2}^{-1}(1-e_{i-1}) \\
	& = \, e_{i-1} + gg_{1}^{-1}\cdots g_{i-2}^{-1}(1-e_{i-1}) 
\end{align*} and hence \begin{align*}
	g_{i} \, &= \, \pi_{e_{i}}\!\left( gg_{1}^{-1}\cdots g_{i-1}^{-1} \right)\! \, = \, e_{i-1} + e_{i}gg_{1}^{-1}\cdots g_{i-2}^{-1}(1-e_{i-1})e_{i} + 1-e_{i} \\
	& = \, e_{i}gg_{1}^{-1}\cdots g_{i-2}^{-1}f_{i} + 1-f_{i} \, = \, f_{i}gg_{1}^{-1}\cdots g_{i-2}^{-1}f_{i} + 1-f_{i} + e_{i-1}gg_{1}^{-1}\cdots g_{i-2}^{-1}f_{i} \\
	&= \, \pi_{f_{i}}\!\left(gg_{1}^{-1}\cdots g_{i-2}^{-1}\right) + e_{i-1}gg_{1}^{-1}\cdots g_{i-2}^{-1}f_{i} \\
	& \stackrel{\ref{lemma:idempotent.difference}}{\in} \, \Gamma_{R}(f_{i}) + e_{i-1}Rf_{i} \, \stackrel{\ref{lemma:small.subgroup}\ref{lemma:small.subgroup.ball}}{\subseteq} \, \trap\!\left(\ball_{1/n}(\GL(R),d_{R})\right)\! .
\end{align*} Since $e_{n} = ({\rk_{R}}\vert_{E})^{-1}(1) = 1$ and thus $\pi_{e_{n}} = \id_{\GL(R_{E})}$, we see that $g_{n} = gg_{1}^{-1}\cdots g_{n-1}^{-1}$. Consequently, \begin{displaymath}
	g \, = \, g_{n}\cdots g_{1} \, \in \, \trap\!\left(\ball_{1/n}(\GL(R),d_{R})\right)^{n} .
\end{displaymath} This proves the first assertion of the theorem.

Since $\GL(R) \cap \A(R)$ is dense in $(\GL(R),d_{R})$ by~\cite[Proposition~8.5]{SchneiderGAFA}, the conclusion above entails that $\GL(R) = \overline{\trap\!\left(\ball_{1/n}(\GL(R),d_{R})\right)^{n}}$. Finally, thanks to~\cite[Lemma~5.4]{BernardSchneider}, \begin{displaymath}
	\ball_{1/2n}(\GL(R),d_{R}) \, \subseteq \, \bigcup \left\{ \Gamma_{R}(e) \left\vert \, e \in \E(R), \, \rk_{R}(e) \leq \tfrac{1}{n} \right\} \! \right. \, \stackrel{\ref{lemma:small.subgroup}\ref{lemma:small.subgroup.ball}}{\subseteq} \, \trap\!\left(\ball_{1/n}(\GL(R),d_{R})\right)\! 
\end{displaymath} and therefore \begin{displaymath}
	\GL(R) \, = \, \overline{\trap\!\left(\ball_{1/n}(\GL(R),d_{R})\right)^{n}} \, \subseteq \, \trap\!\left(\ball_{1/n}(\GL(R),d_{R})\right)^{n+1} . \qedhere
\end{displaymath} \end{proof}

To discuss some ramifications of Theorem~\ref{theorem:bounded}, we now turn to escape functions and the escape property, as introduced in~\cite[Section~3]{SchneiderSolecki24}. Recall that a \emph{length function} on a group $G$ is a map $f \colon G \to \R$ such that \begin{itemize}
	\item[---\,] $f(e) = 0$,
	\item[---\,] $f\!\left( x^{-1}\right) = f(x) \geq 0$ for every $x \in G$, and
	\item[---\,] $f(xy) \leq f(x) + f(y)$ for all $x,y \in G$.
\end{itemize}

\begin{definition}[{\cite[Definition~3.1]{SchneiderSolecki24}}]\label{definition:escape.function} Let $G$ be a topological group and let $\Neigh(G)$ denote the neighborhood filter at the neutral element of $G$. An identity neighborhood $U \in \Neigh(G)$ is said to be an \emph{escape neighborhood} for a length function $f$ on $G$ if \begin{displaymath}
	\forall \epsilon \in \R_{>0} \ \exists n \in \N \colon \quad \{ x \in G \mid x^{1},\ldots,x^{n} \in U \} \, \subseteq \, f^{-1}([0,\epsilon)) .
\end{displaymath} An \emph{escape function} on $G$ is a length function on $G$ possessing an escape neighborhood. We say that $G$ has the \emph{escape property} if, for every $U \in \Neigh(G)$, there exists an escape function $f$ on $G$ such that $f^{-1}([0,1)) \subseteq U$. \end{definition}

\begin{lem}[cf.~{\cite[Proof of Proposition~3.9]{SchneiderSolecki24}}]\label{lemma:escape} Let $G$ be a topological group. If \begin{displaymath}
	\forall U \in \Neigh(G) \ \exists n \in \N \colon \quad G = \overline{\trap(U)^{n}} ,
\end{displaymath} then $G$ does not admit any non-zero escape function. \end{lem}

\begin{proof} Let $f$ be an escape function on $G$ and let $U \in \Neigh(G)$ be an escape neighborhood for $f$. Then $f\vert_{\trap(U)} = 0$. As $f$ is a length function, we conclude that $f\vert_{\trap(U)^{n}} = 0$ for every $n \in \N$. By our hypothesis, there exists $n \in \N$ such that $G = \overline{\trap(U)^{n}}$. Since $f$ is continuous according to~\cite[Remark~3.2]{SchneiderSolecki24}, this implies that $f=0$. \end{proof}

\begin{lem}\label{lemma:homs.vs.escape.functions} Let $G$ and $H$ be topological groups. \begin{enumerate}
	\item\label{lemma:homs.vs.escape.functions.1} If $f$ is an escape function on $H$ and $\pi \colon G \to H$ is a continuous homomorphism, then $f \circ \pi$ is an escape function on $G$.
	\item\label{lemma:homs.vs.escape.functions.2} Suppose that $G$ does not admit any non-zero escape function and that $H$ is Hausdorff and has the escape property. Then every continuous homomorphism from $G$ to $H$ is trivial.
\end{enumerate} \end{lem}

\begin{proof} \ref{lemma:homs.vs.escape.functions.1} This has been observed in~\cite[Item~(a) in the proof of Lemma~3.4]{SchneiderSolecki24}.

\ref{lemma:homs.vs.escape.functions.2} Let $\pi \colon G \to H$ be a continuous homomorphism. If $U$ is an identity neighborhood in $H$, then $H$ admits an escape function $f$ with $f^{-1}([0,1)) \subseteq U$, thus $f \circ \pi$ is an escape function on $G$ by~\ref{lemma:homs.vs.escape.functions.1}, whence $f \circ \pi = 0$ by assumption and so $\pi(G) \subseteq U$. Since $H$ is Hausdorff, this entails triviality of $\pi$. \end{proof}

\begin{cor}\label{corollary:no.escape} If $R$ is a non-discrete irreducible, continuous ring, then $\GL(R)$ does not admit any non-zero escape function. \end{cor}

\begin{proof} This follows by Theorem~\ref{theorem:bounded} and Lemma~\ref{lemma:escape}. \end{proof}

\begin{cor}\label{corollary:homomorphism.rigidity} Let $R$ be a non-discrete irreducible, continuous ring. Every continuous homomorphism from $\GL(R)$ to a Hausdorff topological group with the escape property is trivial. \end{cor}

\begin{proof} This follows from Corollary~\ref{corollary:no.escape} and Lemma~\ref{lemma:homs.vs.escape.functions}\ref{lemma:homs.vs.escape.functions.2}. \end{proof}

We now combine the above with the automatic continuity result of~\cite[Theorem~1.3]{BernardSchneider}. A topological group $G$ is said to have \emph{automatic continuity}~\cite{KechrisRosendal} if every homomorphism from $G$ to any separable topological group is continuous. Moreover, recall that a topological group~$G$ is \emph{$\omega$-narrow}~\cite[Section~2.3, p.~117]{AT} if for every neighborhood $U$ of the neutral element in $G$ there exists a countable subset $C \subseteq G$ such that $UC = G$.

\begin{lem}\label{lemma:automatic.continuity} Let $G$ be a topological group with automatic continuity. Then every homomorphism from $G$ to any $\omega$-narrow topological group is continuous. \end{lem}

\begin{proof} Let $\psi \colon G \to H$ be a homomorphism to an $\omega$-narrow topological group $H$. By work of Guran~\cite{guran} (see also~\cite[Theorem~3.4.23]{AT}), there exists a topological group embedding $\iota \colon H \to \prod_{i \in I} H_{i}$ into the direct product of some family of second-countable (in particular, separable) topological groups $(H_{i})_{i \in I}$. For each $i \in I$, we consider the projection $\pi_{i} \colon \prod_{j \in I} H_{j} \to H_{i}, \, h \mapsto h_{i}$ and note that ${\pi_{i}} \circ \iota \circ \psi$ is continuous thanks to automatic continuity of $G$. Therefore, $\psi$ must be continuous, too. \end{proof}

A group $G$ is called \emph{minimally almost periodic}~\cite{VonNeumann34,HartKunen} (or \emph{weakly mixing}~\cite{BergelsonFurstenberg}) if every homomorphism from $G$ to a compact Hausdorff topological group is trivial. By the work of Peter and Weyl~\cite{PeterWeyl}, any compact Hausdorff topological group embeds into a direct product of finite-dimensional unitary groups (see, e.g.,~\cite[Lemma~1.7]{HartKunen}). Therefore, a group is minimally almost periodic if and only if it has no non-trivial finite-dimensional unitary representation.

\begin{cor}\label{corollary:rigidity} Let $R$ be a non-discrete irreducible, continuous ring. \begin{enumerate}
	\item\label{corollary:rigidity.1} Every action of $\GL(R)$ by isometries on a separable locally compact metric space is trivial. In particular, any action of $\GL(R)$ on a countable set is~trivial.
	\item\label{corollary:rigidity.2} Every homomorphism from $\GL(R)$ to any $\sigma$-compact locally compact topological group is trivial. In particular, $\GL(R)$ is minimally almost periodic.
\end{enumerate} \end{cor}

\begin{proof} \ref{corollary:rigidity.1} Let $X$ be a separable locally compact metric space. Its isometry group $\Iso(X)$, endowed with the topology of pointwise convergence, constitutes a separable Hausdorff topological group possessing the escape property by~\cite[Proposition~3.6(iv)]{SchneiderSolecki24}. Now, any action of $\GL(R)$ by isometries on $X$ gives rise to a homomorphism from $\GL(R)$ to $\Iso(X)$, which then has to be continuous with respect to the rank topology according to~\cite[Theorem~1.3]{BernardSchneider} and thus must be trivial by Corollary~\ref{corollary:homomorphism.rigidity}.

\ref{corollary:rigidity.2} Let $\psi \colon \GL(R) \to G$ be a homomorphism to a $\sigma$-compact locally compact Hausdorff topological group $G$. In particular, $G$ is $\omega$-narrow, thus $\psi$ is continuous thanks to~\cite[Theorem~1.3]{BernardSchneider} and Lemma~\ref{lemma:automatic.continuity}. Since $G$ has the escape property according to~\cite[Proposition~3.6(ii)]{SchneiderSolecki24}, it follows that $\psi$ must be trivial by Corollary~\ref{corollary:homomorphism.rigidity}. (Alternatively, triviality of the continuous homomorphism $\psi$ follows from~\cite[Theorem~2.2.1]{veech} and~\cite[Corollary~1.6]{SchneiderGAFA}, as explained in~\cite[p.~1612, paragraph before Corollary~1.6]{SchneiderGAFA}.) \end{proof}

In order to record some further consequences of Theorem~\ref{theorem:bounded}, let us recall the following piece of terminology. A uniform space $X$ is said to be \emph{bounded} (in the sense of Bourbaki~\cite[II, \S4, Exercise~7, p.~210]{bourbaki}) if, for every entourage $U$ of $X$, there exist $n \in \N$ and a finite subset $F \subseteq X$ such that \begin{displaymath}
	\forall x_{n} \in X \ \exists x_{0} \in F \ \exists x_{1},\ldots,x_{n-1} \in X \ \forall i \in \{ 0,\ldots,n-1\} \colon \quad (x_{i},x_{i+1}) \in U .
\end{displaymath} A uniform space $X$ is bounded if and only if every uniformly continuous real-valued function on $X$ is bounded (see~\cite[Theorem~1.4]{hejcman} or~\cite[Theorem~2.4]{atkin}).

\begin{remark}[{\cite[Theorem~1.20]{hejcman}}]\label{remark:bounded} Every dense subspace of a bounded uniform space is bounded with respect to the relative uniformity. \end{remark}

For any pseudo-metric group $(G,d)$, the uniformity generated by $d$ coincides with both the left and the right uniformity arising from the group topology induced by~$d$. In particular, there is no ambiguity concerning the notion of boundedness for pseudo-metric groups.

\begin{cor}\label{corollary:bounded} Let $R$ be a non-discrete irreducible, continuous ring. Then $\GL(R)$, endowed with the rank topology, is bounded in the sense of Bourbaki. \end{cor}

\begin{proof} Our Theorem~\ref{theorem:bounded} asserts that $\GL(R) = \ball_{1/n}(\GL(R),d_{R})^{n+1}$ for every $n \in \N_{>0}$, which readily entails the desired conclusion. \end{proof}

We conclude this note by discussing a consequence of Corollary~\ref{corollary:bounded} for the unit group of a certain example of a non-discrete irreducible, continuous ring arising naturally from an arbitrary field $K$ as follows. The inductive limit $R$ of the rings \begin{displaymath}
	K \, \cong \, \M_{2^{0}}(K) \, \stackrel{\iota_{0}}{\lhook\joinrel\longrightarrow} \, \ldots \, \stackrel{\iota_{n-1}}{\lhook\joinrel\longrightarrow} \, \M_{2^{n}}(K) \, \stackrel{\iota_{n}}{\lhook\joinrel\longrightarrow} \, \M_{2^{n+1}}(K) \, \stackrel{\iota_{n+1}}{\lhook\joinrel\longrightarrow} \, \ldots
\end{displaymath} relative to the embeddings \begin{displaymath}
	\iota_{n} \colon \, \M_{2^{n}}(K)\,\lhook\joinrel\longrightarrow \,  \M_{2^{n+1}}(K), \quad a\,\longmapsto\, 
	\begin{pmatrix}
		a & 0\\
		0 & a
	\end{pmatrix} \qquad (n \in \N) 
\end{displaymath} is easily seen to be regular. Since ${\rk_{\M_{2^{n+1}}(K)}} \circ {\iota_{n}} = {\rk_{\M_{2^{n}}(K)}}$ for each $n \in \N$, the rank functions of $\M_{2^{n}}(K)$ $(n \in \N)$ jointly extend to a rank function $\rho$ on $R$. The metric completion of $(R,d_{\rho})$, which we denote by $\M_{\infty}(K)$, is a non-discrete irreducible, continuous ring~\cite{NeumannExamples,Halperin68} such that $\cent(\M_{\infty}(K)) = \M_{2^{0}}(K) \cong K$~\cite[Theorem~2.8(c)]{Goodearl78}. Moreover, as established in~\cite[Theorem~1]{Halperin68}, the sequence $(2^{n})_{n \in \N}$ may be replaced by any factor sequence of natural numbers without altering the isomorphism type of the resulting ring $\M_{\infty}(K)$.

\begin{lem}[cf.~{\cite[p.~258]{CarderiThom}}]\label{lemma:density} Let $K$ be a field and let $R \defeq \M_{\infty}(K) = \overline{\bigcup_{n \in \N}\M_{2^{n}}(K)}$. Then $\bigcup\nolimits_{n \in \N} \SL_{2^{n}}(K)$ is dense in $\GL(R)$. \end{lem}

\begin{proof} Let $g \in \GL(R)$ and $\epsilon \in \R_{>0}$. Then we find $n \in \N$ and $a \in \M_{2^{n}}(K)$ such that $2^{-n} \leq \tfrac{\epsilon}{3}$ and $\rk_{R}(g-a) \leq \tfrac{\epsilon}{3}$. Note that $\M_{2^{n}}(K)$ is unit-regular (e.g., due to Remark~\ref{remark:discrete}\ref{remark:discrete.2} and Remark~\ref{remark:unit.regular}). Consequently, thanks to Lemma~\ref{lemma:invertible.approximation}, there exists $b \in \GL(\M_{2^{n}}(K)) = \GL_{2^{n}}(K)$ such that $\rk_{\M_{2^{n}}(K)}(a-b) = 1-\rk_{\M_{2^{n}}(K)}(a)$. Moreover, by Lemma~\ref{lemma:sl.dense}, we find $h \in \SL_{2^{n}}(K)$ with $\rk_{\M_{2^{n}}(K)}(b-h) \leq 2^{-n} \leq \tfrac{\epsilon}{3}$. As ${\rk_{\M_{2^{n}}(K)}} = {{\rk_{R}}\vert_{\M_{2^{n}}(K)}}$ by Remark~\ref{remark:discrete}\ref{remark:discrete.2} and~\cite[Remark~4.5(B)]{BernardSchneider}, we deduce that \begin{align*}
	\rk_{R}(a-b) \, &= \, \rk_{\M_{2^{n}}(K)}(a-b) \, = \, 1-\rk_{\M_{2^{n}}(K)}(a) \, = \, 1-\rk_{R}(a) \\
	& \stackrel{\ref{remark:properties.pseudo.rank.function}\ref{remark:unit.group}}{=} \, \rk_{R}(g)-\rk_{R}(a) \, \stackrel{\ref{lemma:pseudo.rank.function}\ref{lemma:pseudo.rank.function.2}}{\leq} \, \rk_{R}(g-a) \, \leq \, \tfrac{\epsilon}{3}
\end{align*} and thus \begin{displaymath}
	\rk_{R}(g-h) \, \stackrel{\ref{lemma:pseudo.rank.function}\ref{lemma:pseudo.rank.function.2}}{\leq} \, \rk_{R}(g-a) + \rk_{R}(a-b) + \rk_{R}(b-h) \, \leq \, \epsilon .
\end{displaymath} This proves that $\bigcup\nolimits_{n \in \N} \SL_{2^{n}}(K)$ is dense in $\GL(R)$. \end{proof}

A topological group $G$ is called \emph{exotic}~\cite{HererChristensen} if every homomorphism from~$G$ to the unitary group of a Hilbert space, continuous with respect to the strong operator topology, is trivial. A topological group $G$ is said to be \emph{strongly exotic}~\cite{banaszczyk} if every homomorphism from $G$ to the group of invertible bounded linear operators on a Hilbert space, continuous with respect to the weak operator topology, is trivial.

\begin{thm}[\cite{CarderiThom}]\label{theorem:exotic} If $F$ is a finite field, then $\GL(\M_{\infty}(F))$ is exotic. \end{thm}

\begin{cor}\label{corollary:strongly.exotic} If $F$ is a finite field, then $\GL(\M_{\infty}(F))$ is strongly exotic. \end{cor}

\begin{proof} Consider $R \defeq \M_{\infty}(F) = \overline{\bigcup_{n \in \N}\M_{2^{n}}(F)}$. According to Lemma~\ref{lemma:density}, the subgroup $G \defeq \bigcup_{n \in \N} \SL_{2^{n}}(F)$ is dense in $\GL(R)$. In turn, Theorem~\ref{theorem:exotic} entails that $G$ is exotic with respect to the relative topology inherited from $\GL(R)$, as follows from~\cite[Remark~B.2(i)]{SchneiderSolecki24} and~\cite[Proposition~B.1(iii)]{SchneiderSolecki24}. Since $G$ is obviously first-countable with respect to the rank topology, bounded with regard to the same by Corollary~\ref{corollary:bounded} and density in $\GL(R)$, and locally finite and therefore amenable as a discrete group, the topological group $G$ is strongly exotic by~\cite[Lemma~5.2]{SchneiderSolecki24}. Thus, $\GL(R)$ is strongly exotic according to~\cite[Remark~5.1]{SchneiderSolecki24}. \end{proof}

\section*{Acknowledgments}

The author would like to thank Frieder Knüppel for a helpful exchange on the results of Rodgers and Saxl~\cite{RodgersSaxl}. Moreover, the author thanks Andreas Thom for his comments on an earlier version of this manuscript. Finally, the two anonymous referees' valuable remarks and suggestions are gratefully acknowledged.

\end{document}